\documentclass[review,hidelinks,onefignum,onetabnum]{siamart250211}

\usepackage{lipsum}
\usepackage{amsfonts}
\usepackage{epstopdf}
\usepackage{algorithmic,multirow,multicol}
\usepackage{amsmath}
\usepackage{color}
\usepackage{amssymb}
\usepackage{mathrsfs}
\usepackage{latexsym, graphicx,bbm}
\ifpdf
  \DeclareGraphicsExtensions{.eps,.pdf,.png,.jpg}
\else
  \DeclareGraphicsExtensions{.eps}
\fi

\newcommand{\e}{\varepsilon}
\def\mR{{\mathbb R}}

\newsiamremark{remark}{Remark}
\newsiamremark{assumption}{Assumption}

\headers{A class of multi-fidelity methods for the Schr\"odinger equation}{Y. Lin, and L. Liu}

\title{On a class of multi-fidelity methods for the semiclassical Schr\"odinger equation with uncertainties\thanks{Submitted to the editors May 22, 2024.
\funding{L. Liu acknowledges the support by National Key R\&D Program of China (2021YFA1001200), Ministry of Science and Technology in China, Early Career Scheme (24301021) and General Research Fund (14303022 \& 14301423) funded by Research Grants Council of Hong Kong from 2021-2023. Y. Lin was partially supported by National Key R\&D Program of China (2024YFA1016100) and National Natural Science Foundation of China (12201404).
	 }}}

\author{Yiwen Lin\thanks{School of Mathematical Sciences and Institute of Natural Sciences, Shanghai Jiao Tong University, Shanghai 200240, P.R. China
  (\email{linyiwen@sjtu.edu.cn}).}
\and Liu Liu\thanks{The Chinese University of Hong Kong, Hong Kong
  (\email{lliu@math.cuhk.edu.hk}).}}

\usepackage{amsopn}

\ifpdf
\hypersetup{
  pdftitle={A class of multi-fidelity methods for the Schr\"odinger equation},
  pdfauthor={Y. Lin, and L. Liu}
}
\fi

\begin{document}

\maketitle

\begin{abstract}
In this paper, we study the semiclassical Schr\"odinger equation with random parameters and develop several robust multi-fidelity methods. We employ the time-splitting Fourier pseudospectral (TSFP) method for the high-fidelity solver, and consider different low-fidelity solvers including the meshless method like frozen Gaussian approximation (FGA) and the level set (LS) method for the semiclassical limit of the Schr\"odinger equation. With a careful choice of the low-fidelity model, we obtain an error estimate for the bi-fidelity method. We conduct numerous numerical experiments and validate the accuracy and efficiency of our proposed multi-fidelity methods, by comparing the performance of a class of bi-fidelity and tri-fidelity approximations. 
\end{abstract}

\begin{keywords}
Schr\"odinger equation, uncertainty quantification, multi-fidelity method, error estimate
\end{keywords}

\begin{MSCcodes}
35J10, 65M70
\end{MSCcodes}

\section{Introduction}

Uncertainty Quantification (UQ) has drawn an increasing attention over the last decades. 
The study of UQ involves quantitative characterization and management of uncertainties in the inputs of models, which is important to assess, validate and improve these models, then to produce more reliable predictions for the outputs. 

In UQ, the non-intrusive stochastic collocation (SC) method has been popularly used. One key challenge to collocation approaches is the computational cost, as it requires repetitive implementations of the deterministic solver, especially for complex problems or high-dimensional random uncertainties, thus makes the simulation difficult or infeasible. There have been many study in this direction, for example \cite{Babuska, Schwab, NXiu, Webster, Xiu07}. Fortunately, there usually exist some approximated, less complicated low-fidelity models in practice. 
Compared to the high-fidelity models, these low-fidelity models usually contain simplified physics or are simulated on a coarser mesh in physical space. 
The goal of multi-fidelity algorithm in the SC framework is to combine the computation efficiency of low-fidelity models with the high accuracy of high-fidelity models, thus can construct an accurate surrogate for the high-fidelity model and produce reliable results, at a significantly reduced simulation cost. We mention that there are alternative multi-fidelity methods. We refer the readers to a nice survey article \cite{Survey2018} that reviews multi-fidelity methods in the contexts of uncertainty propagation, inference and optimization.
Regarding UQ for kinetic and related problems with random inputs, there is a review \cite{L-Review} and a series of works \cite{DP19, DP20,HPW21} that discusses some recent progress on multi-fidelity methods for kinetic equations. 

In this paper, we consider the dimensionless Schr\"odinger equation with random initial data and potential function: 
	\begin{align}\label{Schro-eqn}
		\left\{
		\begin{array}{l}
			\displaystyle i \e\partial_t \psi = -\frac{\e^2}{2} \partial_{xx}\psi + V(x,z)\psi, \quad t>0 , \, x\in \Omega, \\[8pt]
			\displaystyle \psi(t,x,z) = \psi_0(x,z), 
		\end{array}\right.
	\end{align}
where $\Omega=(a,b)\in \mR$ is a bounded domain equipped with periodic boundary conditions. The parameter $\e$, named the semiclassical parameter, describes the microscopic macroscopic scale ratio and is small within the semiclassical regime (i.e. $0<\e\ll1$). 
From $\psi(t,x, z)\in \mathbb{C}$, one can compute the physical quantities including position density and current density:
	$$\begin{aligned}
		& \rho^{\varepsilon}(t,x,z)=\left|\psi^{\varepsilon}(t,x,z)\right|^2,\quad 		
		 J^{\varepsilon}(t,x,z)=\varepsilon \operatorname{Im}\left(\overline{\psi^{\varepsilon}(t,x,z)} \nabla \psi^{\varepsilon}(t,x,z)\right),
	\end{aligned}$$
	where ``--" denotes complex conjugation.

By combining an inaccurate but cheap low-fidelity solver and accurate yet expensive solvers in clever ways, multi-fidelity methods have shown effectiveness in reducing the number of high-fidelity samples. 
For the high-fidelity solver, we utilize the time-splitting Fourier pseudospectral (TSFP)  method \cite{BaoCai2014, BCF2023} in this paper.  Note that if one uses the stochastic collocation (SC) method for TSFP, the number of collocation points should be $O(1/\e)$ \cite{jin2020gaussian}. It motivates us to develop a class of multi-fidelity method for the Schr\"odinger equation with uncertainties. From numerical tests in Section \ref{sec4.1}, we find that our multi-fidelity method is more efficient than the SC method for TSFP, when achieving the same level of accuracy.
 In order to explore cheaper solvers for the low-fidelity models to significantly reduce the computational cost for the Schr\"odinger equation with random inputs, in this paper we investigate the bi- and tri-fidelity methods for \eqref{Schro-eqn}.

Based on our observations from numerical experiments, we address the motivation of the tri-fidelity method in the perspectives of efficiency and accuracy, through the offline-online two-stage implementation shown in Algorithm 2.1. The tri-fidelity method employs three models: a cheap low-fidelity model $\boldsymbol{U}^{(1)}$, a moderately expensive but more accurate medium-fidelity model $\boldsymbol{U}^{(2)}$, and the expensive high-fidelity model $\boldsymbol{U}^{(H)}$. Here we compare against two bi-fidelity approaches: BF-1 that uses $\boldsymbol{U}^{(2)}$ and $\boldsymbol{U}^{(H)}$, where $\boldsymbol{U}^{(2)}$ serves as the low-fidelity model, and BF-2 that uses $\boldsymbol{U}^{(1)}$ and $\boldsymbol{U}^{(H)}$, with $\boldsymbol{U}^{(1)}$ being the low-fidelity model. These two different choices of low-fidelity model in the bi-fidelity model will help us illustrate how the tri-fidelity method can achieve a superior balance between efficiency and accuracy. 
\begin{itemize}
    \item [i)] The key gain in efficiency comes from the offline point selection stage, where we exclusively use the cheapest model $\boldsymbol{U}^{(1)}$ to select important parameter points in the random space. Compared to the bi-fidelity solver BF-1 that uses  $\boldsymbol{U}^{(2)}$ for point selection, our tri-fidelity approach employs $\boldsymbol{U}^{(1)}$ thus significantly reduce the computational cost, as simulating $\boldsymbol{U}^{(1)}$ requires far less computational effort than $\boldsymbol{U}^{(2)}$.
\item [ii)] In terms of the accuracy improvement, the online stage of coefficient inference utilizes the medium-fidelity model $\boldsymbol{U}^{(2)}$, which is closer to the high-fidelity solution than $\boldsymbol{U}^{(1)}$. Compared to the bi-fidelity method BF-2 that uses $\boldsymbol{U}^{(1)}$ for coefficient inference, this tri-fidelity approach yields more accurate result since $\boldsymbol{U}^{(2)}$ better captures the behavior of the high-fidelity solution in the parameter space.
\end{itemize}
In summary, we strategically allocate the computational resources among different hierarchies of models in the design of tri-fidelity algorithm. 
In particular, we choose the cheapest model $\boldsymbol{U}^{(1)}$ in the offline  stage for important point selection and a more accurate $\boldsymbol{U}^{(2)}$ model in the online coefficient inference stage, in order to achieve an almost optimal result through balancing the trade-offs between efficiency and accuracy of the method. 
This explains the motivation of designing the tri-fidelity method.

Besides efficiency and accuracy, the curse of dimensionality of mesh-based methods and the high oscillation of wave function coming from $t$, $x$ and random variable $z$ are another challenges for the UQ problem of semiclassical Schr\"odinger equation. Therefore, meshless methods such as the frozen Gaussian approximation is studied in this work and we try to seek accurate approximations for the high-fidelity solution by employing a tri-fidelity method.

Regarding the analysis, we give an error estimate in Section \ref{sec:estimate} on the bi-fidelity method for solving the uncertain Schr\"odinger equation \eqref{Schro-eqn}. Nevertheless, in practice, it is hard to predict the performance of multi-fidelity methods, and whether the tri-fidelity methods work better than bi-fidelity method in specific problems and settings. In order to better access and evaluate the quality of the multi-fidelity solution, computing empirical error bounds \cite{GZW20, Hampton2018} becomes a necessary task and we will explore it in the numerical tests. 

This paper is organized as follows. Section \ref{sec:multi} gives an introduction to several solvers for the Schr\"odinger equation as well as its semiclassical limit. Section \ref{sec:estimate} is dedicated to establishing an error estimate of the bi-fidelity method, by considering FGA and TSFP scheme as the low- and high-fidelity solvers, respectively. In Section \ref{sec:numerical}, we compare bi-fidelity methods with different choices of low-fidelity solvers and present numerical results for both bi- and tri-fidelity methods.   Finally, conclusions are shown in Section \ref{sec:conclusion}.

\section{Multi-fidelity method}
    \label{sec:multi}

In many practical problems, there may exist several models that describe similar physics but own different simulation cost. The multi-fidelity approximation \cite{Narayan2014,ZNX} has recently been developed. In this paper, we mainly discuss the bi-fidelity and tri-fidelity approaches. Below we first give a summary of the bi-fidelity algorithm \cite{ZNX}. The sought bi-fidelity approximation of $u^H$ is constructed by \eqref{uB-eqn}. 
In the offline stage, we use the cheap low-fidelity model to select the most important parameter points $\gamma=\{z_1, \cdots, z_r\}$ via the greedy procedure \cite{DeVore}. 
In the online stage, for any sample point $z\in I_z$, we project the low-fidelity solution $\boldsymbol{U}^L(z)$ onto the low-fidelity approximation space $\mathscr{U}^L(\gamma)$ by \eqref{Online-eq}, with projection coefficients $c_k$ computed in \eqref{Gsystem}.

For the tri-fidelity approximation, as proposed in \cite{ZNX}, the key of the tri-fidelity method is to (i) use the low-fidelity model to select the ``important" points in the random space; (ii) then infer the projection coefficients by employing the medium-fidelity solution; (iii) finally construct the tri-fidelity solution by using the same approximation rule as that for the medium-fidelity model. We briefly review the bi- and tri-fidelity algorithm in \cref{alg:2}. Based on the observations obtained from our numerical experiments, we address the motivation of developing tri-fidelity method over bi-fidelity approximations in the introduction above.

\begin{algorithm}[t]
	\caption{Bi-fidelity and tri-fidelity approximation for a high-fidelity solution at given $z$.}	\label{alg:2}
	\begin{algorithmic}
				\STATE{\textbf{Offline:}}
				\STATE{\textbf{1} Select a sample set $\Gamma_N=\left\{z_1, z_2, \ldots, z_N\right\} \subset I_z$.}
				\STATE{\textbf{2} 
                Conduct simulation for the low-fidelity model $\boldsymbol{U}^L(z_j)$ at each  $z_j \in \Gamma_N$ to obtain the low-fidelity solution space  $\mathscr{U}^L(\Gamma_N)$.} \STATE{\textbf{3} 
                Select $K$ ``important" points from $\Gamma_N$ and denote $\gamma_K=\left\{z_{i_1}, \cdots, z_{i_K}\right\} \subset \Gamma_N$. For Bi-fidelity approximation, denote the low-fidelity approximation space by $\mathscr{U}^{\mathcal{Y}}\left(\gamma_K\right)= \text{span}
            \{\boldsymbol{U}^\mathcal{Y}(z_{i_1}), \cdots, \boldsymbol{U}^\mathcal{Y}(z_{i_K})\}$, $\mathcal{Y} = L$; for tri-fidelity approximation, run the medium-fidelity model at each sample point of $\gamma_K$ to obtain $\mathscr{U}^{\mathcal{Y}}\left(\gamma_K\right)= \text{span}
            \{\boldsymbol{U}^\mathcal{Y}(z_{i_1}), \cdots, \boldsymbol{U}^\mathcal{Y}(z_{i_K})\}$, $\mathcal{Y} = M$. }
				\STATE{\textbf{4} 
                Run the high-fidelity model at each sample point of $\gamma_K$.
                Denote the high-fidelity approximation space by 
                $\mathscr{U}^H\left(\gamma_K\right)=\text{span} 
            \{\boldsymbol{U}^H(z_{i_1}), \cdots, \boldsymbol{U}^H(z_{i_K})\}$.}
				\STATE{\textbf{Online:}}
				\STATE{\textbf{5} 
                For any given $z \in I_Z$, compute the low-fidelity ($\mathcal{Y} = L$) or medium-fidelity ($\mathcal{Y} = M$)  solution $\boldsymbol{U}^{\mathcal{Y}}(z)$, and the corresponding projection coefficients: 
			\vspace{-0.6em}
				\begin{align}\label{Online-eq}
					\boldsymbol{U}^{\mathcal{Y}}(z)\approx \mathcal{P}_{\mathscr{U}^{\mathcal{Y}}\left(\gamma_K\right)}
                \left[\boldsymbol{U}^{\mathcal{Y}}(z)\right]= \sum_{k=1}^K c_k^{\mathcal{Y}}(z) \boldsymbol{U}^{\mathcal{Y}}\left(z_{i_k}\right),  \ \mathcal{Y} = \{L, M\},
			\end{align}				
                \vspace{-0.6em}
                where $\mathcal{P}_{\mathscr{U}^{\mathcal{Y}}\left(\gamma_K\right)}$ is a projection operator onto the space $U^{\mathcal{Y}}\left(\gamma_K\right)$ with projection coefficients $\left\{c_k^{\mathcal{Y}}\right\}$ computed by the Galerkin approach:
			\vspace{-0.6em}
				\begin{align}\label{Gsystem}
				    \mathbf{G}^{\mathcal{Y}} \mathbf{c}^{\mathcal{Y}}=\mathbf{g}^{\mathcal{Y}},  \quad \mathbf{g}^{\mathcal{Y}}=\left(g_k^{\mathcal{Y}}\right)_{1 \leq k \leq K}, \quad g_k^{\mathcal{Y}}=\left\langle \boldsymbol{U}^{\mathcal{Y}}(z), \boldsymbol{U}^{\mathcal{Y}}\left(z_{i_k}\right)\right\rangle,
				\end{align}
                \vspace{-0.6em}
                where $\mathbf{G}^{\mathcal{Y}}$ is the Gramian matrix for $\mathscr{U}^{\mathcal{Y}}(\gamma_K)$,
                \begin{align}\label{Gramian}
                    (\mathbf{G}^{\mathcal{Y}})_{lk} = \left\langle \boldsymbol{U}^{\mathcal{Y}}\left(z_{i_{\ell}}\right), \boldsymbol{U}^{\mathcal{Y}}\left(z_{i_k}\right)\right\rangle^{\mathcal{Y}},  \ 1 \leq \ell, k \leq K,
                \end{align}
				and $\langle\cdot, \cdot\rangle^{\mathcal{Y}}$ is the inner product associated with $\mathscr{U}^{\mathcal{Y}}\left(\gamma_K\right)$.}
				\STATE{\textbf{6} 
                Construct the tri-fidelity approximation by applying the same approximation rule as in the low-(medium-) fidelity model: 
        		\vspace{-0.6em}
				\begin{equation}\label{uB-eqn}
					\boldsymbol{U}^{\mathcal{F}}(z)=\sum_{k=1}^K c_k^{\mathcal{Y}}(z) \boldsymbol{U}^H\left(z_{i_k}\right), \  \mathcal{F} = \{B, T\}, \  \mathcal{Y} = 
                    \{L, M\}.
			\vspace{-0.6em}
		\end{equation}}
	\end{algorithmic}
\end{algorithm}

    In this paper, we will always employ the time-splitting Fourier pseudospectral (TSFP) method as the high-fidelity method, the frozen Gaussian approximation (FGA) as one choice of the low-fidelity solvers and level set (LS) method as another choice of the low-fidelity solvers. For the rest of this section, we will give a brief introduction of these solvers.

\subsection{The high-fidelity TSFP method}
\label{sec:tsfp}

In this section, we will introduce our choice of the high-fidelity solver for solving the Schr\"odinger equation \eqref{Schro-eqn}, namely the time-splitting Fourier pseudospectral (TSFP) method \cite{BCF2023}. First we mention some existing numerical schemes for the deterministic Schr\"odinger equation. 
 
Due to the highly oscillatory structure of the wave function, various effective numerical algorithms have been proposed and developed. Traditional algorithms can be divided into two classes: direct discretization methods and asymptotic methods \cite{engquist2003, jin2008gaussian, zhou2014numerical, jin2020gaussian,miao2023novel}.
Popular choices of direct discretization methods include finite difference time domain method, the time-splitting spectral method (TSSP) \cite{BJM2002} and time-splitting Fourier pseudospectral (TSFP) method \cite{BaoCai2014, BCF2023}. 

Define the spatial mesh size $h=\Delta x>0$ with $h=(b-a) / M$ for $M$ an even positive integer and choose the time step $\tau=\Delta t>0$. Then let the grid points and the time step be
	$
	x_j:=a+j h,  t_n:=n \tau,  j=0,1, \ldots, M,  n=0,1,2, \ldots.
	$	
	For any given $z$, let $\Psi_j^{\varepsilon, n}$ be the approximation of $\psi^{\varepsilon}\left(x_j, t_n,z\right)$. From time $t=t_n$ to time $t=t_{n+1}$, the Schr\"odinger equation \eqref{Schro-eqn} can be decomposed into the following two subproblems via the Strang splitting: 
	\begin{align}
		i \e\partial_t \psi + \frac{\e^2}{2} \partial_{xx}\psi = 0, \label{TSFPstep1} \\
		i \e\partial_t \psi - V(x,z)\psi = 0.\label{TSFPstep2}
	\end{align}
	Equation \eqref{TSFPstep1} will be discretized in space by the spectral method and integrated in time exactly.  Equation \eqref{TSFPstep2} will then be solved exactly. 
The detailed TSFP method for discretizing the Schr\"odinger equation \eqref{Schro-eqn} can be given for $n \geq 0$ as \cite{BCF2023}
\begin{align}\label{TSFP}
	\begin{aligned}
		& \Psi_j^{(*)}=\sum_{l \in \mathcal{T}_N} e^{-i \frac{\varepsilon\tau \mu_l^2}{4}} (\widetilde{\Psi^n})_l e^{i \mu_l\left(x_j-a\right)}, \\
		& \Psi_j^{(**)}=e^{-i \frac{\tau}{\varepsilon} V\left(x_j,z\right)} \Psi_j^{(*)}, \quad\quad\quad\quad\quad\quad\quad j = 0,1,2,\ldots,N, \\
		& \Psi_j^{n+1}=\sum_{l \in \mathcal{T}_N} e^{-i \frac{\varepsilon\tau \mu_l^2}{4}} (\widetilde{\Psi^{(**)}})_l e^{i \mu_l\left(x_j-a\right)},
	\end{aligned}
\end{align}
where 
$\widetilde{\Psi}_l^n=\frac{1}{N} \sum_{j=0}^{N-1} \Psi_j^n e^{-i \mu_l\left(x_j-a\right)}, \quad \mu_l=\frac{2 \pi l}{b-a}, \ l =-\frac{N}{2}, \ldots, \frac{N}{2}-1,
$
with $\Psi_j^0=\Psi_0\left(x_j\right)$ for $j = 0,1,2,\ldots,N$.
Let $I_N u$ be the trigonometric interpolation of a function $u$ \cite{Shen2011spectral}, 
$$I_N u=\sum_{l=-N/2}^{N/2-1} \widetilde{u}_l e^{i \mu_l(x-a)}, \  \widetilde{u}_l=\frac{1}{N} \sum_{j=0}^{N-1} u(x_j) e^{-i \mu_l\left(x_j-a\right)}, \  l =-\frac{N}{2}, \ldots, \frac{N}{2}-1.
$$	
As a remark, we only focus on the case with one-dimensional spacial variable $(d=1)$. Generalization to $d>1$ is straightforward for tensor product grids and the results remain valid without modifications. 

\subsection{Low-fidelity solvers}

It is known that the low-fidelity models are not unique, in this section we will present several of our choices for the low-fidelity models. For convenience, we review the following low-fidelity solvers for the corresponding deterministic problems. 

\subsubsection{Frozen Gaussian Appproximation (FGA) method}
    \label{sec:fga}

We first mention the advantages of the Frozen Gaussian Appproximation (FGA) or Frozen Gaussian Sampling (FGS) methods. 
The Eulerian-based methods suffer from the ``curse of dimensionality", thus making them highly computationally expensive \cite{jin2011mathematical,lasser2020computing}. For example, in the time-splitting spectral method, the degree of freedom to approximate the wave function on a $d$-dimensional mesh is $O(\e^{-d})$. 
On the other hand, the so-called FGA \cite{Lu2017frozen} or FGS methods \cite{FGS3,FGS1,FGS2} provide us favorable schemes for solving high-dimensional spatial problems. 
		Denote by $z_0=\left(q_0, p_0\right)$ the phase space variables. Assume that $A\left(0, z_0\right)$ is an integrable function in $\mathbb{R}^{2 d}$, i.e.,
		$
		\int_{\mathbb{R}^{2 d}} \left|A\left(0, z_0\right)\right| \mathrm{d} z_0<\infty.
		$
		The usual FGA on a single surface \cite{Lu2017frozen} is given by
		$$	\begin{aligned}
			\psi_{\mathrm{FGA}}(t, x) & =\frac{1}{(2 \pi \varepsilon)^{3 d / 2}} \int_{\mathbb{R}^{2 d}}  A\left(t, z_0\right) \exp \left(\frac{i}{\varepsilon} \Theta\left(t, x, z_0 \right)\right) \mathrm{d} z_0\\
			& =\frac{1}{(2 \pi \varepsilon)^{3 d / 2}} \int_{\mathbb{R}^{2 d}} \left|A\left(0, z_0\right)\right| \frac{A\left(t, z_0\right)}{\left|A\left(0, z_0\right)\right|} \exp \left(\frac{i}{\varepsilon} \Theta\left(t, x, z_0\right)\right)\mathrm{d} z_0,
		\end{aligned} $$
		and the initial condition
		$$
		\psi_0(x)=\frac{1}{(2 \pi \varepsilon)^{3 d / 2}} \int_{\mathbb{R}^{2 d}}  A\left(0, z_0\right) \exp \left(\frac{i}{\varepsilon} \Theta\left(0, x,  z_0\right)\right) \mathrm{d} z_0,
		$$
	where
	$$
	\begin{aligned}
		\Theta(t, x, q, p) & =S(t, q, p)+P(t, q, p) \cdot(x-Q(t, q, p))+\frac{i}{2}|x-Q(t, q, p)|^2, \\
		A(t, q, p) & =a(t, q, p) \int_{\mathbb{R}^d} \psi_0(y) e^{\frac{i}{\varepsilon}\left(-p \cdot(y-q)+\frac{i}{2}|y-q|^2\right)} \mathrm{d} y.
	\end{aligned}
	$$
	Here, Gaussian profile $Q$ and momentum function $P$ are governed by the Hamiltonian flow associated to the classical Hamiltonian $h(q, p)=\frac{1}{2}|p|^2+V(q)$. Action function $S$ is associated to the Hamiltonian flow and the amplitude $a$ is governed by corresponding evolution equations.

Define $T_{n: 1}=\left(t_n, t_{n-1}, \ldots, t_1\right)$ a partition of the time interval $[0, t]$ satisfying
	$0 \leqslant t_1 \leqslant t_2 \leqslant \cdots \leqslant t_n \leqslant t$, with $t_0=0$. When $t \in\left[t_k, t_{k+1}\right)$ for $k$ being an integer, the evolution equations are given accordingly as
		\begin{align}\label{odes}
		    \begin{aligned}
			\frac{\mathrm{d}}{\mathrm{d} t} Q^{(k)} & =P^{(k)}, \\
			\frac{\mathrm{d}}{\mathrm{d} t} P^{(k)} & =-\nabla V(Q^{(k)}), \\
			\frac{\mathrm{d}}{\mathrm{d} t} S^{(k)} & =\frac{1}{2}(P^{(k)})^2-V(Q^{(k)}), \\
			\frac{\mathrm{d}}{\mathrm{d} t} A^{(k)} & =\frac{1}{2} A^{(k)} \operatorname{tr}\left((Z^{(k)})^{-1}\left(\partial_z P^{(k)}-i \partial_z Q^{(k)} \nabla_Q^2 V(Q^{(k)})\right)\right).
		\end{aligned}
		\end{align}
	See Section 4 in \cite{Lu2017frozen} for the asymptotic derivation of these equations.

\subsubsection{The semiclassical limit and Liouville equation}
\label{sec:ls}

When the Wigner transformation is applied to the Schr\"odinger equation, its semi-classical limit reveals a kinetic equation in phase space, which is known as the classical Liouville or Vlasov equation. In the semi-classical regime where the scaled Planck constant $\e$ is small, the Liouville equation is a good approximation for the Schr\"odinger equation, thus one can choose it as a candidate of our low-fidelity models. 
The physical quantities of interests are position and current densities. When the $\e$ is presumed small, the densities computed from the Schr\"odinger equation is an $O(\e)$ approximation to that obtained by the Vlasov equation using moments of the distribution function. By the Wigner transformation \cite{Bal, BJM2002, Wigner,jin2011mathematical}, one can derive that
    for the semiclassical WKB initial data
	$\psi(x, 0)=A_0(x) e^{ \mathrm{i} S_0(x) / \varepsilon}$, the limit of the Wigner transform of $\psi$ given by $w$  satisfies the following Liouville equation:
	\begin{align}\label{liouville-eta}
		\begin{aligned}
			& \partial_t w+p \cdot \nabla_x w-\nabla_x V \cdot \nabla_p w=0, \\
			& w(0, x, p)=\rho_0(x) \delta\left(p-\nabla_x S_0\right) .
		\end{aligned}
	\end{align}

	 The level set (LS) method proposed in \cite{Jin2005JCP} solves the Liouville equation \eqref{liouville-eta} by decomposing $w$ into $f$ and $\phi_i$, namely $w(t, x, p):=f(t, x, p) \Pi_{i=1}^d\delta(\phi_i(t,x,p)),$
     each of which satisfies a Liouville equation with following initial conditions:
	\begin{align}\label{LS}
	    \begin{aligned}
		& \partial_t f+p \cdot \nabla_x f-\nabla_x V(x) \cdot \nabla_p f=0, \\
		& \partial_t \phi_i+p \cdot \nabla_x\phi_i-\nabla_x V(x) \cdot \nabla_p \phi_i=0,\\
		& f(0, x, p)=\rho_0(x) , \ \phi_i(0, x, p)=p-\nabla_x S_0,
	\end{aligned}
	\end{align}
	Then the physical observables of the Liouville equation \eqref{liouville-eta} are thus given by:
$$
\begin{aligned}
	 \rho(t, x)=\int f(t, x, p) \Pi_{i=1}^d \delta\left(\phi_i\right) d p, \quad
	 J(t, x)=\int p f(t, x, p)  \Pi_{i=1}^d \delta\left(\phi_i\right) d p.
\end{aligned}
$$

\section{Error estimates}
	\label{sec:estimate}

        In this section, we will establish an error estimate of a bi-fidelity method for the semiclassical Schr\"odinger equation with random parameters. 
 
	Let the high-fidelity solution be the macroscopic moments of density and momentum that are obtained from the solution $\Psi$ to the Schr\"odinger equation \eqref{Schro-eqn} solved via TSFP method.	
	Then define the high-fidelity solution $\mathbf{U}^H$ of the Schr\"odinger equtaion \eqref{Schro-eqn} by the solution $\Psi^H$, i.e.,
	$$ \mathbf{U}^H = |\Psi^H|^2 := \mathbf{\rho}^H. $$
	Let the low-fidelity solution be the macroscopic moments obtained from the solution $\Psi$ to the Schr\"odinger equation \eqref{Schro-eqn} solved via FGA method.	
	Then define the low-fidelity solution $\mathbf{U}^L$ of the Schr\"odinger equation \eqref{Schro-eqn} by the solution $\Psi^L$, i.e.,
	$$ \mathbf{U}^L = |\Psi^L|^2 := \mathbf{\rho}^L. $$

	Let the bi-fidelity solution $\mathbf{U}^B$ \eqref{uB-eqn} solved from the bi-fidelity approximation Algorithm 1.	
	In order to obtain the error estimate of $\boldsymbol{U}^H-\boldsymbol{U}^B$ in general, we use the following way to split the total error, by inserting the information of $u^L$ :
	\begin{align}\label{uHB}
			\begin{aligned}
			& \boldsymbol{U}^H(z)-\boldsymbol{U}^B(z) 
			=  \boldsymbol{U}^H(z)-\sum_{k=1}^K c_k(z) \boldsymbol{U}^H\left(z_k\right) \\
			= & \boldsymbol{U}^H(z)-\boldsymbol{U}^L(z)+\left(\boldsymbol{U}^L(z)-\sum_{k=1}^K c_k(z) \boldsymbol{U}^L\left(z_k\right)\right)+\sum_{k=1}^K c_k(z)\left(\boldsymbol{U}^L\left(z_k\right)-\boldsymbol{U}^H\left(z_k\right)\right),
		\end{aligned}
	\end{align}
	where the second term is nothing but the projection error of the greedy algorithm, and it remains to estimate $\boldsymbol{U}^H(z)-\boldsymbol{U}^L(z)$ in proper norms.

	\subsection{Notations and prerequisites}
    We first define the space and norms that will be used. The Hilbert space of the random variable is given by
	$$
	H\left(\mathbb{R}^d ; \pi \mathrm{d} z\right)=\left\{f \mid I_z \rightarrow \mathbb{R}, \int_{I_z} f^2(z) \pi(z) \mathrm{d} z<\infty\right\},
	$$
	and equipped with the inner product
	$
	\langle f, g\rangle_\pi=\int_{I_z} f g \pi(z) \mathrm{d} z.
	$
	We introduce the standard multivariate notation. 
	Denote the countable set of "finitely supported" sequences of non-negative integers by
	$$
	\mathcal{F}:=\left\{\nu=\left(\nu_1, \nu_2, \cdots\right): \nu_j \in \mathbb{N} \text {, and } \nu_j \neq 0 \text { for only a finite number of } j\right\}
	$$
	with $|\nu|:=\sum_{j \geq 1}\left|\nu_j\right|$. For $\nu \in \mathcal{F}$ supported in $\{1, \cdots, J\}$, the partial derivative in $z$ is defined by
	$
	\partial_z^\nu=\frac{\partial^{|\nu|}}{\partial^{\nu_1} z_1 \cdots \partial^{\nu_J} z_J} .
	$
	The $z$-derivative of order $\nu$ of a function $f$ is denoted by $f^\nu=\partial_z^\nu f$.

\bigskip

Now we introduce some prerequisites that will be used later for the FGA and TSFP numerical schemes. In the following two lemmas, for simplicity of notations, we consider the deterministic problem setting. The error estimate for FGA is reviewed in \cref{thm:FGA} under \cref{assump:A}; the error estimate for TSFP is obtained based on \cite{BCF2023} and summarized in \cref{thm:L2TSFP} under \cref{AssumpAB}. We give the proof of \cref{thm:L2TSFP} in \cref{app:A}. 
 
	\begin{assumption}\label{assump:A}
		 $V(q) \in C^{\infty}\left(\mathbb{R}^d\right)$ satisfies $
	\sup _{q \in \mathbb{R}^d}\left|\partial_\alpha V(q)\right| \leq C_V, \  \text{for  } |\alpha|=2 .
	$
	\end{assumption}
Assume that there exists a compact set $K \subset \mathbb{R}^{2 d}$ such that the initial approximation error, given by
\begin{align}\label{def-eps-in}
	\epsilon_{\text {in }}=\left\|\frac{1}{(2 \pi \varepsilon)^{3 d / 2}} \int_K A^{(0)}(0, z) e^{\frac{i}{\varepsilon} \Phi^{(0)}(0, x, z)} \mathrm{d} z-u_0(0, x)\right\|_{L^2\left(\mathbb{R}^d\right)} > 0, 
\end{align}
is sufficiently small for the accuracy requirement. 

 \begin{lemma} \label{thm:FGA}
		\cite{Lu2017frozen} Let $\Psi_{FGA}(t, x)$ be the approximation given by the FGA (with phase space integral restricted to the set $K$) for the Schr\"odinger equation \eqref{Schro-eqn}, whose exact solution is denoted by $\psi(t, x)$. Under \cref{assump:A}, for any given final time $t$, there exists a constant $C$, such that for any $\varepsilon>0$ sufficiently small, we have
		$$
		\left\|\Psi_{F G A}(t, x)-\psi(t, x)\right\|_{L^2\left(\mathbb{R}^d\right)} \leqslant C \varepsilon+\epsilon_{\text {in }},
		$$
		where $\epsilon_{\text {in }}$ is the initial approximation error defined in \eqref{def-eps-in}.
	\end{lemma}

	\begin{assumption}\label{AssumpAB}
  Assume that the exact solution $\psi(t,x)$ of the Schr\"odinger equation and the potential satisfy $$\|\psi(t,x)\|_{L^{\infty}\left(\left[0, T_{\varepsilon}\right] ; H_{\mathrm{per}}^m\right)} \leq \frac{C}{\e^m},\quad V(x) \in H_{\mathrm{per}}^{m^*}, \  m^*=\max \{m, 4\},$$
  where $m$ describes the regularity of the exact solution. Here, $H_{\mathrm{per }}^m(\Omega)=\{\psi \in \left.H^m(\Omega) \mid \partial_x^k \psi(a)=\partial_x^k \psi(b), \ k=0,1, \ldots, m-1\right\},$ with the equivalent $H^m$-norm on $H_{\mathrm{per}}^m(\Omega)$ given as $\|\psi\|_{H^m}=(\sum_{l \in \mathbb{Z}}(1+\mu_l^2)^m|\widehat{\psi}_l|^2)^{1 / 2}$. 
	\end{assumption}	
	\begin{lemma}\label{thm:L2TSFP}
		Let $\psi^n$ be the numerical approximation obtained from the TSFP \eqref{TSFP}. Under \cref{AssumpAB}, for any $0<\varepsilon \leq 1$, $t_n\in[0,T]$, we have
		\begin{align}\label{H1TSFP}
			\left\|\psi\left(t_n, x\right)-I_N \psi^n\right\|_{L^2} \leq C_0 G_m \frac{T}{\varepsilon} \left(h^{m}+\frac{\tau^2}{\varepsilon^3}\right), m\geq 2,
		\end{align}
		 where $C_0$ and $C_1$ are positive constants independent of $\varepsilon, h, \tau, n, m$; $G_m$ is a positive constant independent of $\varepsilon, h, \tau, n$ and only depends on $\|\psi\|_{L^{\infty}\left([0, T] ; L^2\right)}$ and $\|V\|_{L^2}$.
	\end{lemma}

\subsection{Main result}

    Now we introduce the main result on the bi-fidelity method in \cref{thm:bi-fi} under \cref{assump:A},  \cref{AssumpAB} and \cref{Assumph}. Here TSFP is used as the high-fidelity solver and FGA is chosen as the low-fidelity model. 
    
	\begin{assumption}\label{Assumph}
		Assume each component of the random variable $z:=(z_j)_{j \geq 1}$ has a compact support. Let $\left(\psi_j\right)_{j \geq 1}$ be an affine representer of the random initial data $\psi_{\text {in }}$, which by definition means that \cite{Cohen2015}
		\begin{align}\label{Assump:hin}
			\psi_{\text {in }}(z)=\tilde{\psi}_0+\sum_{j \geq 1} z_j \phi_j,
		\end{align}
		where $\tilde{\psi}_0=\tilde{\psi}_0(x)$ is independent of $z$, and the sequence $\left\{\left\|\phi_j\right\|_{L^{\infty}(\Omega)}\right\}_{j \geq 1} \in \ell^p$ for $0<p<1$, with $\Omega$ representing the physical space.
	\end{assumption}	
	The reason we need to assume a compact support of random variable $z$ and \eqref{Assump:hin} is because of \cite{Cohen2015}, which is used for the projection error estimate for greedy algorithm. We now show the main result:	
	\begin{theorem}\label{thm:bi-fi}
		Let \cref{assump:A},  \cref{AssumpAB},  \cref{Assumph} hold, and 
  $\boldsymbol{U}^H$, $\boldsymbol{U}^B$ represent the high- and bi-fidelity approximation for the position density solved from the Schr\"odinger equation \eqref{Schro-eqn}. 
		With $K$ high-fidelity simulation runs, the error estimate of the bi-fidelity method is given by
		$$
		\left\|\boldsymbol{U}^H(z)-\boldsymbol{U}^B(z)\right\|_{L_x^2 L_z^2} \leq \frac{C_1}{\e} \frac{1}{(K/2+1)^{q/2}}+ C_2 \sqrt{\frac{K}{\lambda_0}}\,  \left( \varepsilon +\frac{T}{\varepsilon}\left(h^m+\frac{\tau^2}{\varepsilon^3}\right)+ \epsilon_{in}\right),
		$$
		where $C_1, C_2$ are constants independent of $\varepsilon$; $m$ depends on the regularity of the exact solution $\psi$ as shown in \cref{AssumpAB}; $\lambda_0$ is the minimum eigenvalue of the Gramian matrix defined in \eqref{Gramian}; $q$ is independent of dimensionality of $z$;  $h$ and $\tau$ are spatial and temporal sizes.
	\end{theorem}

Now we show the proof for \cref{thm:bi-fi}.

\textbf{Proof.}
	Taking $||\cdot||_{L_x^2 L_z^2}$ norm on both sides of the equality \eqref{uHB} for the moments, one has
	\begin{align}\label{uHBnorm}
			\begin{aligned}
				& \left\|\boldsymbol{U}^H(z)-\boldsymbol{U}^{B}(z)\right\|_{L_x^2 L_z^2}=\left\|\boldsymbol{U}^H(z)-\sum_{k=1}^K c_k^{L}(z) \boldsymbol{U}^H\left(z_k\right)\right\|_{L_x^2 L_z^2} \\
				& \leq \underbrace{\left\|\boldsymbol{U}^H(z)-\boldsymbol{U}^L(z)\right\|_{L_x^2 L_z^2}}_{\operatorname{Term} A_1}+\underbrace{\left\|\boldsymbol{U}^L(z)-\sum_{k=1}^K c_k^{L}(z) \boldsymbol{U}^L\left(z_k\right)\right\|_{L_x^2 L_z^2}}_{\text {Term } A_2} \\
				& +\underbrace{\left\|\sum_{k=1}^K c_k^{L}(z)\left(\boldsymbol{U}^L\left(z_k\right)-\boldsymbol{U}^H\left(z_k\right)\right)\right\|_{L_x^2 L_z^2}}_{\text {Term } A_3} .
			\end{aligned}
	\end{align}
 Let $\boldsymbol{U}^{H,L}:=\boldsymbol{\rho}^{H,L}$ be the high- or low-fidelity approximation for the position density, and assume $\mathbf{u}=\rho$ be the exact solution obtained from $\psi$ that solves the Schr\"odinger equation. The first term $A_1$ is estimated by $\left\|\boldsymbol{U}^H-\mathbf{u}\right\|_{L_x^2 L_z^2}$ and $\left\|\boldsymbol{U}^L-\mathbf{u}\right\|_{L_x^2 L_z^2}$. The second term $A_2$ is actually the projection error of the greedy algorithm when searching the most important points $\gamma_K$ from the low-fidelity solution manifold. To estimate the third term $A_3$, a bound for the vector $\|\mathbf{c}\|$ with $\|\cdot\|$ the matrix induced $\ell_2$ norm is obtained and then the Cauchy-Schwarz inequality is applied.
	
Now that
$$ \left\|\mathbf{U}^H(z)-\mathbf{U}^L(z)\right\|_{L_x^2 L_z^2}
	\leq \left\|\mathbf{U}^H-\mathbf{u}\right\|_{L_x^2 L_z^2} + \left\|\mathbf{u}-\boldsymbol{U}^L\right\|_{L_x^2 L_z^2}. $$
Since $\boldsymbol{U} = \boldsymbol{\rho} = |\psi|^2$, then 
$||\boldsymbol{U}||_{L_x^2 L_z^2} = ||\boldsymbol{\rho}||_{L_x^2 L_z^2} \leq ||\psi||_{L_x^2 L_z^2}^2 $. 
Let $\Psi_{TSFP}(t,x,z)$ be the numerical approximation obtained from the TSFP. Let $\Psi_{FGA}(t, x,z)$ be the approximation solved by the FGA  and $\psi(t, x,z)$ be the exact solution. Then, the following estimate holds from \cref{thm:L2TSFP} and \cref{thm:FGA}:
\begin{align}\label{TermA1}
    \begin{aligned}
		||\mathbf{U}^L_{FGA} -\mathbf{U}^H_{TSFP} ||_{L_x^2 L_z^2}
		&\leq ||\mathbf{U}^L_{FGA} - \mathbf{u} ||_{L_x^2 L_z^2} + || \mathbf{u}- \mathbf{U}^H_{TSFP} ||_{L_x^2 L_z^2} \\
		& \lesssim 
		  || \Psi_{FGA}^L - \psi ||_{L^2_x L_z^2}  +  || \psi - \Psi_{TSFP}^H ||_{L^2_x L_z^2} \\
		&\leq  C  \left( \varepsilon +\frac{T}{\varepsilon}\left(h^m+\frac{\tau^2}{\varepsilon^3}\right)\right) + \epsilon_{in}. 
	\end{aligned}
\end{align}
The second term $A_2$ in \eqref{uHBnorm} is evaluated by the Kolmogorov $K$-width of a functional manifold,  which characterizes the optimal distance for approximation from a general $K$-dimensional subspace \cite{GambaLiu2021, JinLin2023bifidelity}. Denote by $d_K(\boldsymbol{U}^L(I_z))$ the Kolmogorov $K$-width of the functional manifold $\boldsymbol{U}^L\left(I_z\right)$, defined by
$$
d_K(\boldsymbol{U}^L(I_z))=\inf _{\operatorname{dim}\left(V_K\right)=K} \sup _{v \in \boldsymbol{U}^L\left(I_z\right)} d^L\left(v, V_K\right).$$ 
Denote the space $\mathcal{H}=L_x^2$. 
Specifically, one has
\begin{align}
	\begin{aligned}		d_K(\boldsymbol{U}^L(I_z))_{\mathcal{H}} & \leq \sup _{v \in \boldsymbol{U}^L\left(I_z\right)} \min _{w \in V_K}\|v-w\|_{\mathcal{H}}=\sup _{z \in I_z} \min _{w \in V_K}\left\|\boldsymbol{U}^L(z)-w\right\|_{\mathcal{H}} \\
		& \leq\left\|\boldsymbol{U}^L-\sum_{\nu \in \Lambda_K} w_\nu P_\nu\right\|_{L^{\infty}\left(I_z, \mathcal{H}\right)} \leq \frac{C}{\varepsilon}\frac{1}{(K+1)^q}, \quad q=\frac{1}{p}-1,
	\end{aligned}
\end{align}
which holds for all $z$. Then
\begin{align}\label{dK-estimate}
   \left\|\boldsymbol{U}^L-\sum_{\nu \in \Lambda_K} w_\nu P_\nu\right\|_{L_z^2\mathcal{H}} \leq \frac{C}{\varepsilon}\frac{1}{(K+1)^q}, \quad q=\frac{1}{p}-1,
\end{align}
based on the analysis of the bound for $\partial_z\rho$ in Appendix C.
Here $p$ is a constant associated with the affine representer $\left(\psi_j\right)_{j \geq 1}$ in the random initial data specified by Eq. \eqref{Assump:hin}. In addition, $\sum_{\nu \in \Lambda_K} w_\nu P_\nu$ is the truncated Legendre expansion with $\left(P_k\right)_{k \geq 0}$ the sequence of renormalized Legendre polynomials on $[-1,1]$, and $\Lambda_K$ is the set of indices corresponding to the $K$ largest $\left\|w_\nu\right\|_{\mathcal{H}}$. 
Plugging $K / 2$ into Eq. \eqref{dK-estimate},
\begin{align}\label{TermA2}
	\operatorname{Term} A_2 \leq \frac{C_1}{\e} \frac{1}{(K / 2+1)^{q / 2}}.
\end{align}
For the third term $A_3$, one has
\begin{align}\label{ck-norm}
	\left(\sum_{k=1}^K\left\|c_k(z)\right\|_{L_z^2}^2\right)^{1 / 2} \lesssim \frac{1}{\sqrt{\lambda_0}}\left\|\boldsymbol{U}^L(z)\right\|_{L_x^2 L_z^2},
\end{align}
where $\lambda_0$ is the smallest eigenvalue of the Gramian matrix $\mathbf{G}_L$ defined in \eqref{Gramian}.
Details for \eqref{ck-norm} are given in \cref{app:B}.
By the Cauchy-Schwarz inequality, combining \eqref{TermA1} and \eqref{ck-norm}, one gets 
\begin{align}\label{TermA3}
	\begin{aligned}
		\operatorname{Term} A_3 & \leq\left(\sum_{k=1}^K\left\|c_k(z)\right\|^2\right)^{1 / 2}\left(\sum_{k=1}^K\left\|\boldsymbol{U}^H\left(z_k\right)-\boldsymbol{U}^L\left(z_k\right)\right\|_{L_x^2 L_z^2}^2\right)^{1/2} \\
		& \lesssim\left(\sum_{k=1}^K\left\|c_k(z)\right\|_{L_z^2}^2\right)^{1/2}\left(\sum_{k=1}^K\left\|\boldsymbol{U}^H\left(z_k\right)-\boldsymbol{U}^L\left(z_k\right)\right\|_{L_x^2}^2\right)^{1/2} \\
		& \leq \sqrt{N}\left(\sum_{k=1}^K\left\|c_k(z)\right\|_{L_z^2}^2\right)^{1/2}\left(\max_k\left\|\boldsymbol{U}^H\left(z_k\right)-\boldsymbol{U}^L\left(z_k\right)\right\|_{L_x^2}\right) \\
		& \leq \frac{\sqrt{K}}{\sqrt{\lambda_0}}\, \left(C'  \left( \varepsilon +\frac{T}{\varepsilon}\left(h^m+\frac{\tau^2}{\varepsilon^3}\right)\right) + \epsilon_{in}\right) \left\|\boldsymbol{U}^L(z)\right\|_{L_x^2 L_z^2} \\
  & \leq C \sqrt{\frac{K}{\lambda_0}}\,  \left( \varepsilon +\frac{T}{\varepsilon}\left(h^m+\frac{\tau^2}{\varepsilon^3}\right)+ \epsilon_{in}\right).
	\end{aligned}
\end{align}
Finally, by adding up \eqref{TermA1}, \eqref{TermA2} and \eqref{TermA3}, the proof is completed from \eqref{uHBnorm}.


\begin{remark}
 Numerically we found that the bi-fidelity error is smaller than the low-fidelity error, this is because the last two terms in equation (3.1) may work together to counteract the first term of low-fidelity error due to the projection effect and thus make the total bi-fidelity error smaller. We would like to emphasize that even though our error bound may not be sharp, making the bound smaller than the low-fidelity error is {\it not} the motivation of this error analysis. 
 {\it The goal of our error bound analysis is to provide some theoretical guarantee that the bi-fidelity error is bounded by some small value.} The core motivation of the multi-fidelity method is to reduce the computational cost significantly and tackle the curse of dimensionality in UQ problems. Compared to the traditional method such as Monte-Carlo or standard stochastic collocation method, in the multi-fidelity algorithm we utilize the computationally cheaper low-fidelity solver to select a few points in random space, then one only needs to run the high-fidelity solver on these few points, combined with the negligible cost of computing the project coefficients. As demonstrated in numerical experiments, for example in Test 1 one can see that by $O(10)$ runs of the high-fidelity solver, we can obtain a  satisfactory level of accuracy at $\mathcal{O}(10^{-5})$. 
\end{remark}

\begin{remark}
Note that if we want to get an error estimate for the tri-fidelity method, then
\begin{align}\label{uHTnorm}
			\begin{aligned}
				& \left\|\boldsymbol{U}^H(z)-\boldsymbol{U}^{T}(z)\right\|_{L_x^2 L_z^2}=\left\|\boldsymbol{U}^H(z)-\sum_{k=1}^K c_k^{M}(z) \boldsymbol{U}^H\left(z_k\right)\right\|_{L_x^2 L_z^2} \\
				& \leq \underbrace{\left\|\boldsymbol{U}^H(z)-\boldsymbol{U}^M(z)\right\|_{L_x^2 L_z^2}}_{\operatorname{Term} A_1}+\underbrace{\left\|\boldsymbol{U}^M(z)-\sum_{k=1}^K c_k^{M}(z) \boldsymbol{U}^M\left(z_k\right)\right\|_{L_x^2 L_z^2}}_{\text {Term } A_2} \\
				& +\underbrace{\left\|\sum_{k=1}^K c_k^{M}(z)\left(\boldsymbol{U}^M\left(z_k\right)-\boldsymbol{U}^H\left(z_k\right)\right)\right\|_{L_x^2 L_z^2}}_{\text {Term } A_3} .
			\end{aligned}
	\end{align}
 Compared to \eqref{uHBnorm}, we substitute $\boldsymbol{U}^{L}$ by $\boldsymbol{U}^{M}$ in \eqref{uHTnorm}. Thus if {\small  $\left\|\boldsymbol{U}^H(z)-\boldsymbol{U}^{M}(z)\right\|_{L_x^2 L_z^2}$ $\leq \left\|\boldsymbol{U}^H(z)-\boldsymbol{U}^{L}(z)\right\|_{L_x^2 L_z^2}$}, then 
    {\small $\left\|\boldsymbol{U}^H(z)-\boldsymbol{U}^{T}(z)\right\|_{L_x^2 L_z^2}\leq \left\|\boldsymbol{U}^H(z)-\boldsymbol{U}^{B}(z)\right\|_{L_x^2 L_z^2}$}, which indicates that the tri-fidelity method is more accurate than the bi-fidelity method if the same low-fidelity solver is chosen in both the bi- and tri-fidelity methods. Moreover, the ``K" important points $\{z_{i_1}, z_{i_2},\ldots, z_{i_K}\}$ in the tri-fidelity method are chosen based on the low-fidelity solver, which is computationally much cheaper than the medium-fidelity solver. Therefore it is necessary to develop multi-fidelity methods.
\end{remark}

\section{Numerical examples}
\label{sec:numerical}

To examine the performance and accuracy of our proposed methods, the numerical errors are defined below. We choose a fixed set of points $\left\{\hat{z}_i\right\}_{i=1}^N \subset I_z$ that is independent of the point sets $\Gamma$, and evaluate the following error between the bi-(or tri-) fidelity and high-fidelity solutions at a final time $t$ :
\begin{eqnarray}\label{def:err}
	\mathcal{E} \approx \frac{1}{N} \sum_{i=1}^N\left\{\dfrac{1}{N_x}\sqrt{\sum_{j=1}^{N_x}\left|u^H\left(t, x_j, \hat{z}_i\right)-u^M\left(t, x_j, \hat{z}_i\right)\right|^2}\right\},  
\end{eqnarray}
where $u^M$ represents the bi-fidelity or tri-fidelity approximations, and $\|\cdot\|_{L^2(\mathcal{X})}$ is the $L^2$ norm in the physical domain $\mathcal{X}$. The error can be considered as an approximation to the average $L^2$ error in the whole space of $\mathcal{X} \times I_z$.

Since our goal is to numerically compute the Schr\"odinger equation with random inputs, thus the high-fidelity model is always chosen as the TSFP method for solving the Schr\"odinger equation  described in Section \ref{sec:tsfp}. The low-fidelity solver will be considered as the FGA method discussed in Section \ref{sec:fga} or the LS method stated in Section \ref{sec:ls}.

	The CFL condition for the Liouville equation is
	$
	\displaystyle\Delta t \leqslant \frac{\Delta x}{2 \max _{(x, p) \in \Omega}\left(p, \nabla_x V\right)}.
	$
	If we employ the LS method in our bi (tri)-fidelity solvers, we will use the simple piecewise linear kernel
	$$
	\begin{aligned}
		& \delta_\eta^{(1)}(x)= \begin{cases}\frac{1}{\eta}\left(1-\frac{|x|}{\eta}\right), & \left|\frac{x}{\eta}\right| \leqslant 1, \\[4pt]
			0, & \left|\frac{x}{\eta}\right|>1,\end{cases} 
   \end{aligned}\  
   \begin{aligned}
		& \delta_\eta^{(2)}(x)= \begin{cases}\frac{1}{2 \eta}\left(1+\cos \frac{\mid \pi x \mid}{\eta}\right), & \left|\frac{x}{\eta}\right| \leqslant 1, \\[4pt]
			0, & \left|\frac{x}{\eta}\right|>1,\end{cases}
	\end{aligned}
	$$
	which have the so-called ``exact integration property": 
	\begin{align*}
	    \sum_{j=-N}^N \delta_{\kappa_0 h}^{(k)}\left(x_j-x_0\right) h=1, \quad \text { for any }-N h<x_0<N h, \quad k=1,2 .
	\end{align*}
	In our numerical tests, we use the kernel $\delta_\eta^{(2)}(x)$ and approximate the gradients $\nabla_x \phi$, $\nabla_p \phi$ by the upwind scheme with a fifth order WENO approximation, For time discretization, the 3rd order TVD Runge-Kutta scheme \cite{BJM2002} is adopted. For the FGA solver, the time evolution ODEs \eqref{odes} is computed by the forth-order Runge-Kutta scheme \cite{Lu2017frozen}. 
 
Without loss of generality, the $d$-dimensional random variables $\mathbf{z}$ are assumed to follow the uniform distribution on $[-1,1]^d$. Let the size of training set $\Gamma$ to be $M=1000$, 
and we evaluate the errors of bi (tri)-fidelity approximation with respect to the number of high-fidelity runs by computing the norm defined in \eqref{def:err} over an independent set of $N=1000$ Monte Carlo samples. 

 In our test problems, the initial condition in \eqref{Schro-eqn} is always chosen as the classical WKB form
$$
u^{\varepsilon}(t=0,x)=u_0^{\varepsilon}(x)=\sqrt{n_0(x)} e^{i S_0(x) / \varepsilon},
$$
with real valued $n_0$ and $S_0$ independent of $\varepsilon$, regular enough, and with $n_0(x)$ decaying to zero sufficiently fast as $|x| \rightarrow \infty$. We choose an appropriate interval $[a, b]$ for the computations such that the periodic boundary conditions do not introduce a significant error relative to the whole space problem.
A Neumann boundary condition is presumed in the $p$ direction. 

Inspired by the empirical error bound estimation approach proposed by \cite{GZW20} to evaluate the performance of the bi- and tri-fidelity approximation,
we give the following empirical error bound estimation approach for the tri-fidelity approximation. The notation corresponds to that in Algorithm \ref{alg:2}.
The advantage of the error bound estimation is that one only needs the low- (or medium-)fidelity data and the first $k+1$ pre-selected high-fidelity samples.
Denoting by $\mathcal{F} = \{B, T\}$ and $\mathcal{Y} = \{L, M\}$, the relative error bound estimation of the bi-(or tri-)fidelity approximation, given $z_*$ can be written as follows:
\begin{align}\label{empirical_bound}
    \frac{\left\|u^H\left(z_*\right)-u^{\mathcal{F}}\left(z_*\right)\right\|}{\left\|u^H\left(z_*\right)\right\|} \leq \frac{d^{\mathcal{Y}}\left(u^{\mathcal{Y}}\left(z_*\right), U^{\mathcal{Y}}\left(\gamma_k\right)\right)}{\left\|u^{\mathcal{Y}}\left(z_*\right)\right\|}\left(c_1+c_2 R_e\left(z_{k+1}\right)\right), 
\end{align}
with
$$
R_e(z):=\frac{\left\|P_{U^H\left(\gamma_k\right)} u^H(z)-u^{\mathcal{F}}(z)\right\|}{d^H\left(z, U^H\left(\gamma_k\right)\right)}.
$$
In our numerical experiments, our empirical results suggest that this error bound estimation is effective, and the constants $c_1$, $c_2$ are both set to be 1.

\subsection{Test 1}
\label{sec4.1}

The initial condition is taken as
{\small$$
\begin{aligned}
	&	n_0(x,\mathbf{z})=\left(\exp \left(-50\left(1+0.8 \sum_{k=1}^{d_1} \frac{z_k^p}{2 k}\right)\left(x-1+0.2 \sum_{k=1}^{d_1} \frac{z_k^q}{2 k}\right)^2\right)\right)^2, \quad\\ &S_0(x,\mathbf{z})=-\frac{1}{5} \ln \left(\exp\left(5(x-1+0.2 \sum_{k=1}^{d_1} \frac{z_k^r}{2 k})\right)+\exp\left(-5(x-1+0.2 \sum_{k=1}^{d_1} \frac{z_k^s}{2 k})\right)\right).
\end{aligned}
$$}
This example was based on \cite{BJM2002}, with uncertainties added to the initial condition. We solve on the $x$-interval $[0,2]$ with periodic boundary conditions. Let $V(x)=10$ be a constant potential. 

Here $\mathbf{z}^R=\left(z_1^R, \cdots, z_{d_1}^R\right)$ with $R = \{p,q,r,s\}$ represents the random variables shown in the initial condition. We set $d_1=5$, then the total dimension $d$ of the random space is 20 in our problem settings. 
In Example 1, we vary $\varepsilon$ from $1/64$ to $1/256$ and compare two different low-fidelity solvers. 
We employ the FGA method or LS method as the low-fidelity solver in Case I and II, respectively. In both cases, the TSFP method is chosen as the high-fidelity solver. 
When $\varepsilon=1/64$, we utilize the FGA method with spatial mesh size $\Delta x_{FGA} = 0.0312$ for the low-fidelity solver in Case I; and the LS method with $\Delta x_{LS} = 0.01, \Delta p = 0.1$ for the low-fidelity solver in Case II. The range of $p$ in this example is $[-2,2]$. 
In TSFP method, $\Delta x_{TSFP} = 0.001$ is chosen. 
When $\varepsilon = 1/256$, we let the low-fidelity solver 1 (FGA) with $\Delta x_{FGA} = 0.0078$, the low-fidelity solver 2 (LS) with $\Delta x_{LS} = 0.0026, \Delta p = 0.1$, and the high-fidelity model (TSFP) with $\Delta x_{TSFP} = 0.00026$. 
We take $\Delta t = 10^{-5}$ in all these methods, 500 number of particles when $\varepsilon=1/64$ and 1300 number of particles when $\varepsilon=1/256$ in the FGA solver.

\begin{figure}[!t]
	\centering
	\includegraphics[width=0.38\textwidth]{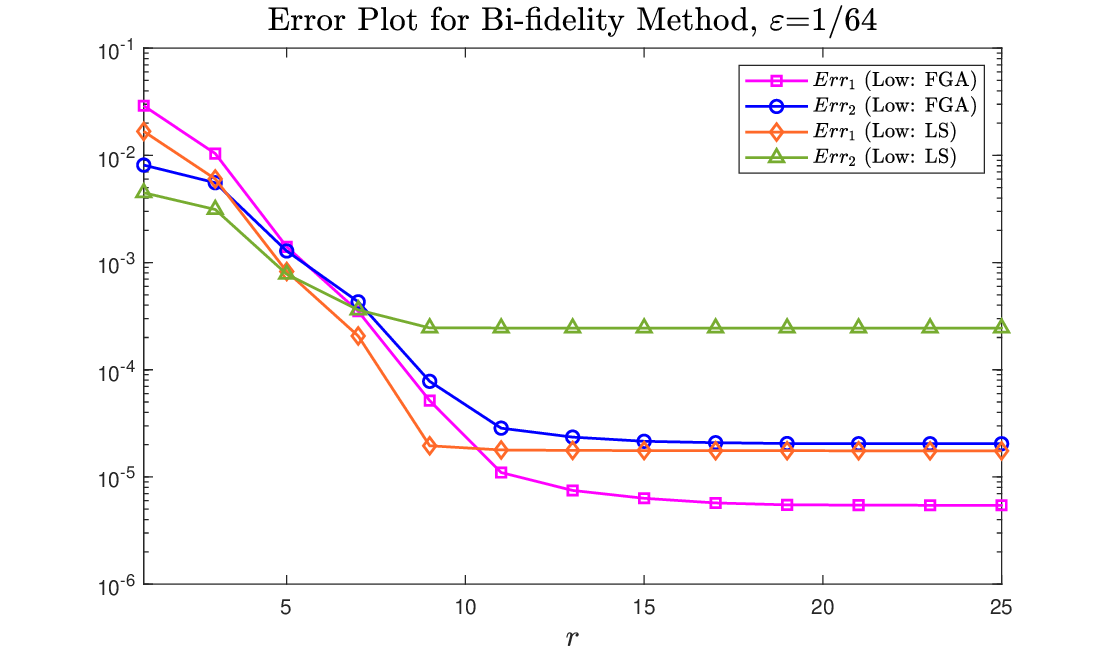}
	\includegraphics[width=0.38\textwidth]{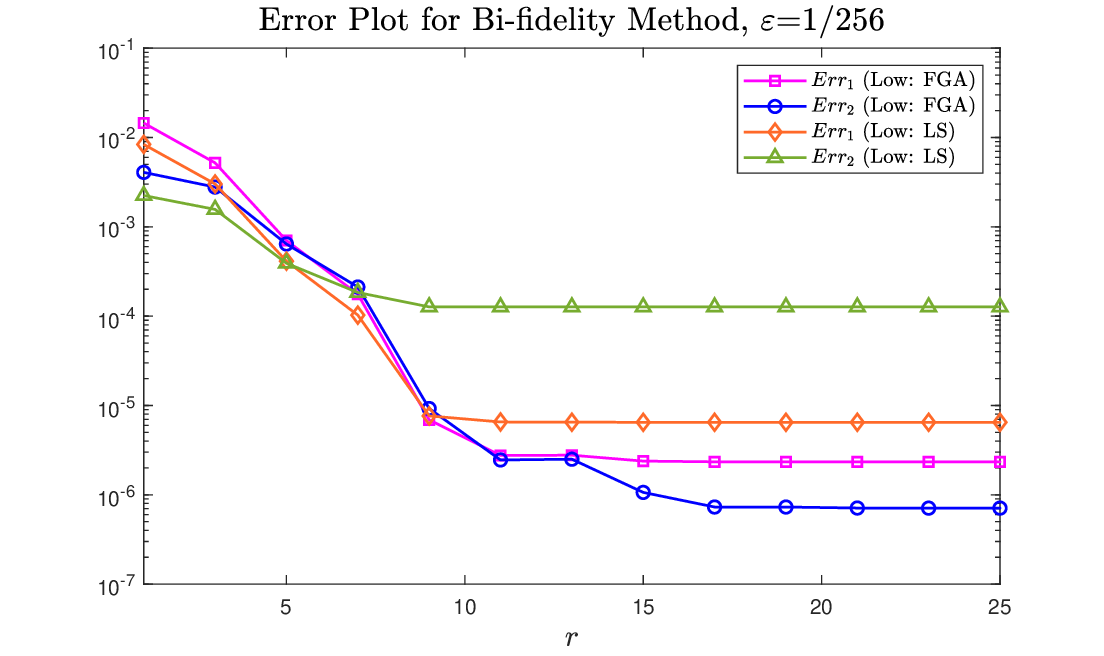}
	\caption{Test 1: Case I and Case II. The mean $L^2$ error of the bi-fidelity approximation $\rho, J$ with respect to the number of high-fidelity runs for different $\varepsilon$. }\label{Ex2Fig1Err}
\end{figure}

First, we investigate the performance of the bi-fidelity approximation for two kinds of low-fidelity solvers. From Figure \ref{Ex2Fig1Err}, fast convergence of the mean $L^2$ errors between the high- and bi-fidelity solutions with respect to the number of high-fidelity runs is observed. Here $\text{Err}_1$ in the figures below stands for the errors of $\rho$, and $\text{Err}_2$ stands for the errors of $J$. This is corresponding to the error estimate analysis in Section \ref{sec:estimate}, where we denote the bi-fidelity solution by 
$\mathbf{u} = (\rho, J)^T $. 
With only 10 high-fidelity runs, the bi-fidelity approximation can reach an accuracy level of  $\mathcal{O}\left(10^{-3}\right)$ for a 20-dimensional problem in random space, which is quite satisfactory compared to other sparse grid type of methods.
We can observe that the smaller $\varepsilon$ is, the lower level the errors ($Err_1, Err_2$) achieve when they saturate. Here $Err_1$, $Err_2$ stand for the errors calculated by \eqref{def:err} for $\rho$ and $J$ respectively. For example, comparing the two orange lines in the left and right part of Figure \ref{Ex2Fig1Err}, the $Err_1$ decreases when $\e$ becomes smaller, since the discrepancy between the low-fidelity model chosen as the Liouville equation and the high-fidelity model known by the Schr\"odinger equation is reduced when $\e$ becomes smaller.

To illustrate the efficiency of the bi-fidelity method, we compare with the stochastic collocation (SC) method \cite{Xiu2005} here. Let $\left\{\mathbf{z}^{(m)}\right\}_{m=1}^{N_c} \subset I_{\mathbf{z}}$ be the set of collocation nodes, $N_c$ the number of samples in the random space. For each individual sample $\mathbf{z}^{(m)}$, one applies the TSFP solver to the deterministic Schr\"odinger equation and obtains the solution ensemble $\textbf{u}\left(t, x, \mathbf{z}^{(m)}\right)$, then adopts the Legendre-Gauss quadrature rule to construct the approximation for the mean value $\bar{\textbf{u}}_{N_c}(t,x) = \sum_{m=1}^{N_c} \textbf{u}\left(t,x,(\mathbf{z}^{(m)}\right)) w_j$, where $w_j$ are the corresponding quadrature weights.   
We let the dimensionality of $z$ to be $d_1=1$ and $N_{ref} = 256$ as a reference solution. We plot the error 
$$Err_{j,N_c} = ||\bar{\textbf{u}}_{j,N_c} - \bar{\textbf{u}}_{j,N_{ref}} ||_{l^2}, \quad j= 1,2,$$ ranging from $N_c = 8$ to $N_c = 128$ as shown in Table \ref{table1}. One can conclude that more than a hundred collocation points are needed to achieve $10^{-4}$ accuracy level for the SC method, while we just need to run the high-fidelity solver in the online stage of $O(10)$ times, in order to achieve the same level of accuracy. 

\begin{table}[t]\label{table1}
\centering
\begin{tabular}{l|l|l|l|l|l|l}
\hline
& $N_c$          & 8      & 16     & 32     & 64       & 128      \\ \hline
\multirow{2}{*}{$\e=1/64$}                        
& $Err_1$  & 0.1583 & 0.0418 & 0.0106 & 0.0026 & 5.17e-04 \\ \cline{2-7} 
& $Err_2$  & 0.0970 & 0.0257 & 0.0066 & 0.0016 & 3.19e-04 \\ \hline
\multicolumn{1}{c|}{\multirow{2}{*}{$\e=1/256$}} 
& $Err_1$ & 0.3607 & 0.0951 & 0.0242 & 0.0058 & 0.0012 \\ \cline{2-7} 
\multicolumn{1}{c|}{}                           
& $Err_2$ & 0.2455 & 0.0650 & 0.0166 & 0.0040 & 8.06e-04 \\ \hline
\end{tabular}
\caption{\text{Tab. 1.} Relative Error of the stochastic Collocation method for TSFP.}
\end{table}

\begin{figure}[!t]
	\centering
	\includegraphics[width=0.49\textwidth]{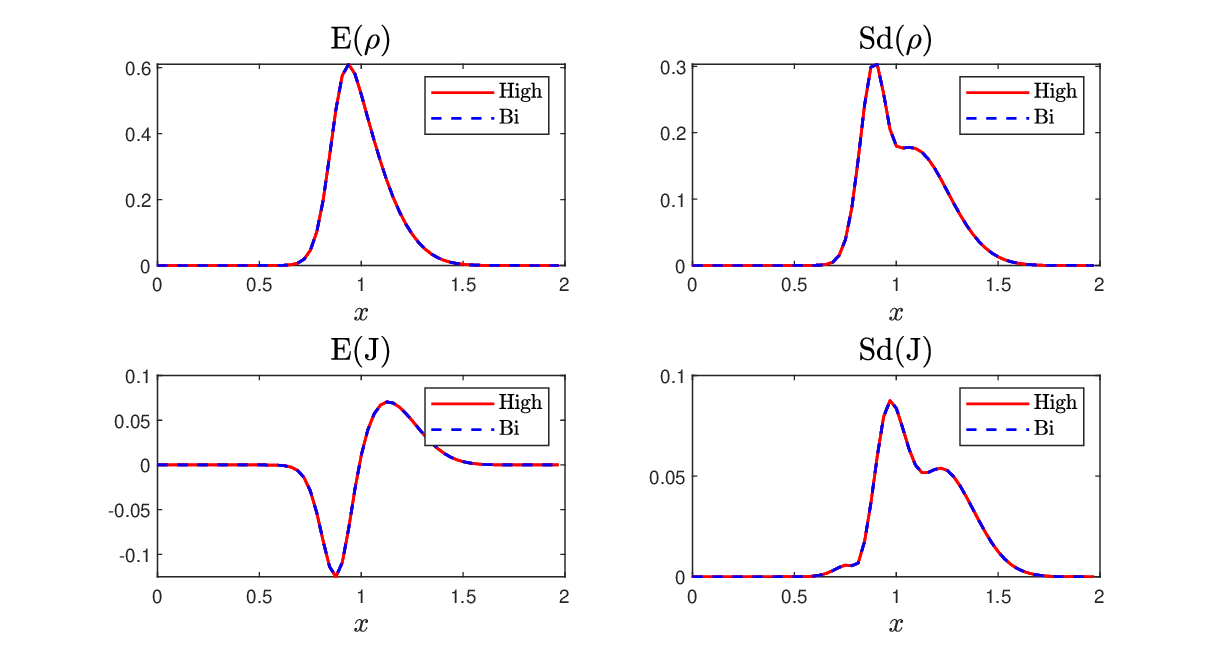}
	\caption{Test 1: Case I. The mean and standard deviation of high- and bi- fidelity solutions $\rho$ and $J$ with $r=25, \varepsilon=1/256$ at $t=0.5s$ with FGA and TSFP for the low- and high-fidelity solver.   }\label{Ex2Fig2N256}
\end{figure}

\begin{figure}[!t]
	\centering 
 \includegraphics[width=0.49\textwidth]{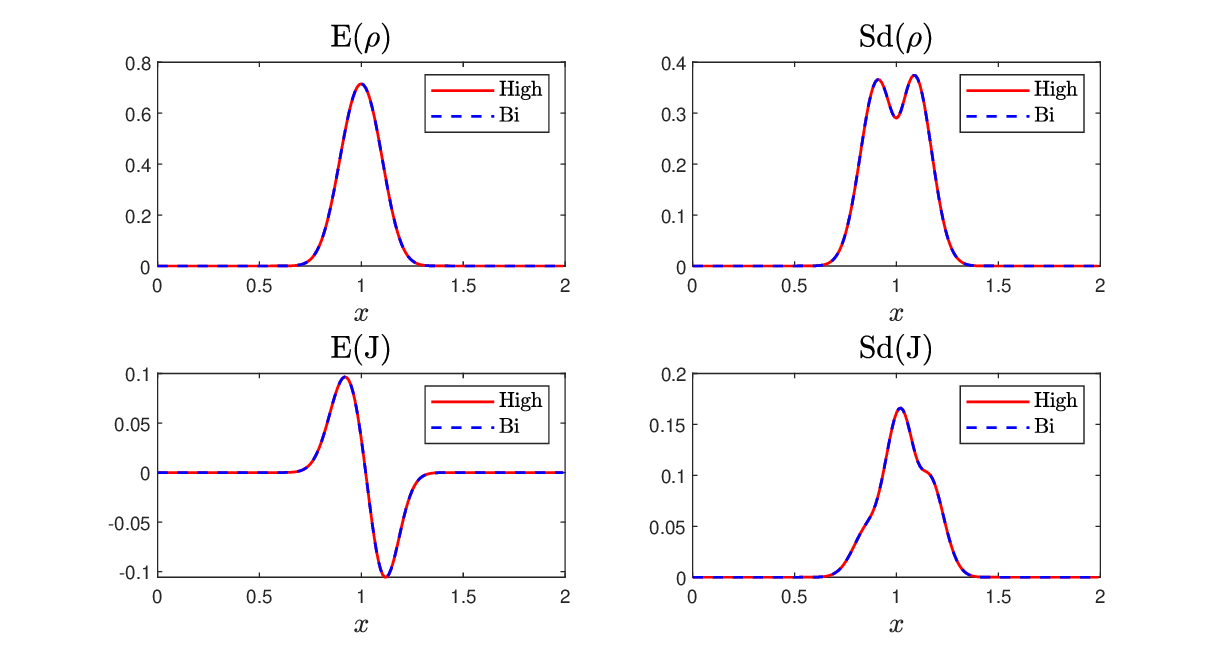}
	\caption{Test 1: Case II. The mean and standard deviation of high- and bi-fidelity solutions $\rho$ and $J$ with $r=25, \varepsilon=1/64$ at $t=0.01s$ with LS and TSFP for the low- and high-fidelity solver.}\label{Ex2Fig2N64}
\end{figure}

Figure \ref{Ex2Fig2N256} and Figure \ref{Ex2Fig2N64}  show clearly that the mean and standard deviation of the bi-fidelity approximation for $\rho$ and $J$ agree well with the high-fidelity solutions by using only 20 high-fidelity runs for both two cases (FGA as low-fidelity for Case I and LS as low-fidelity for Case II).  
The result is a bit surprising yet reasonable, suggesting that although coarser mesh grid in the physical space is taken in both cases, the bi-fidelity approximation can still capture important variations of the high-fidelity model in the random space.

\begin{figure}[!t]
	\centering
	\includegraphics[width=0.49\textwidth]{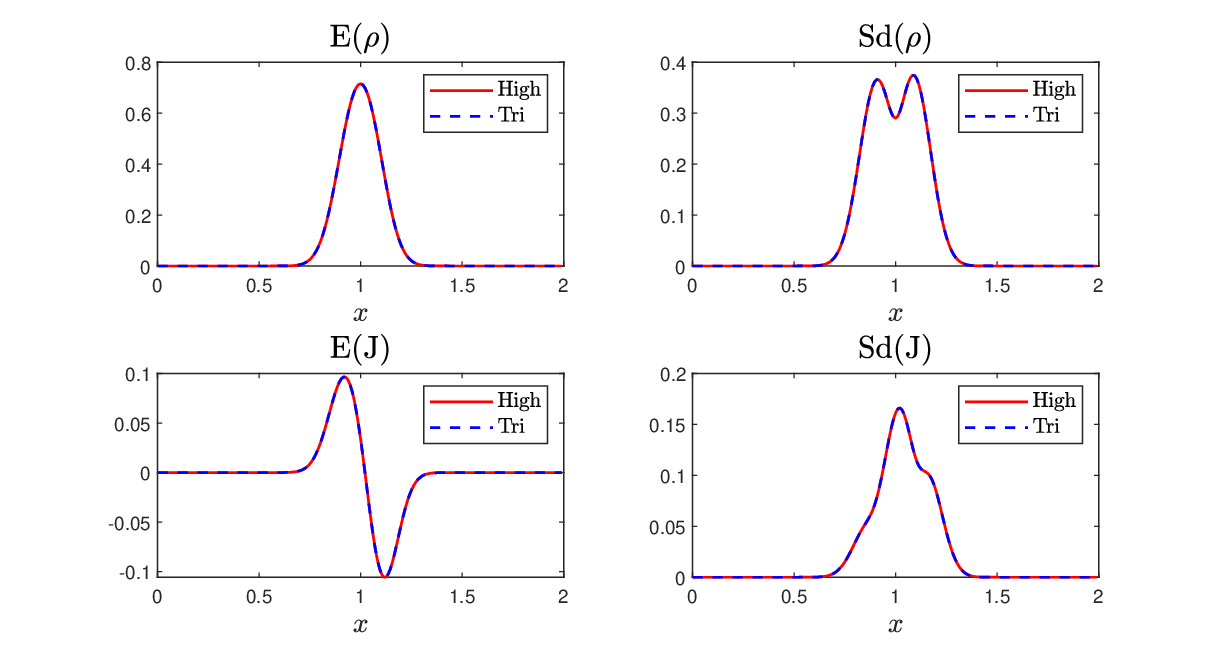}
	\includegraphics[width=0.49\textwidth]{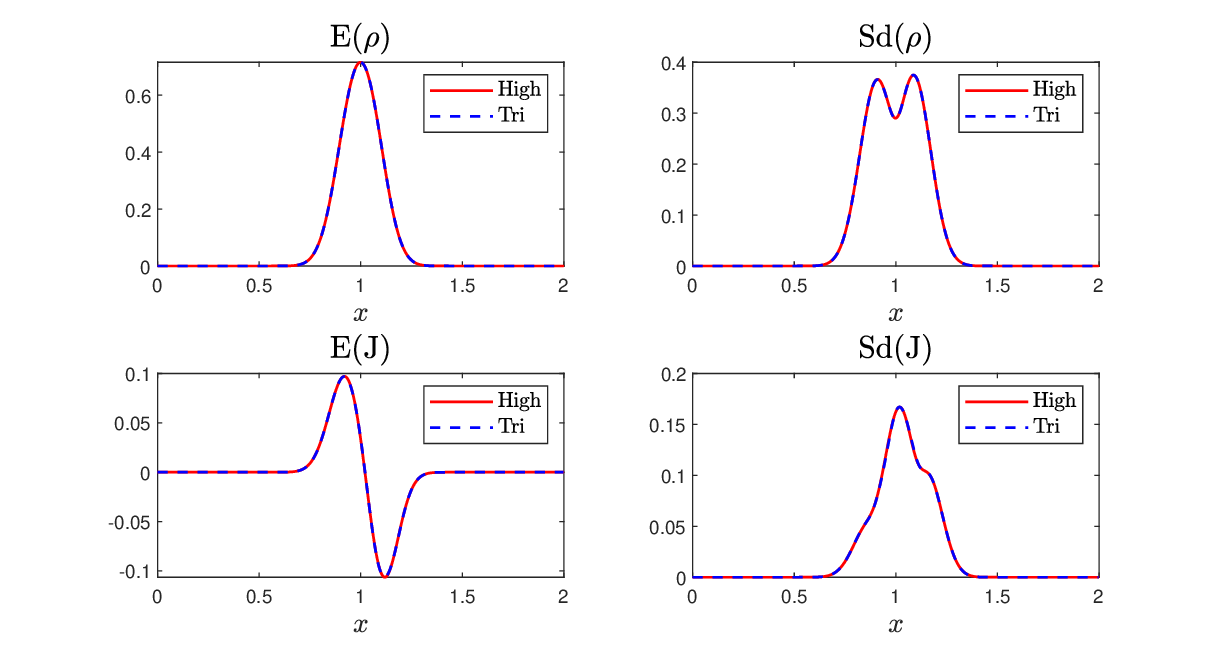}
    \includegraphics[width=0.49\textwidth]{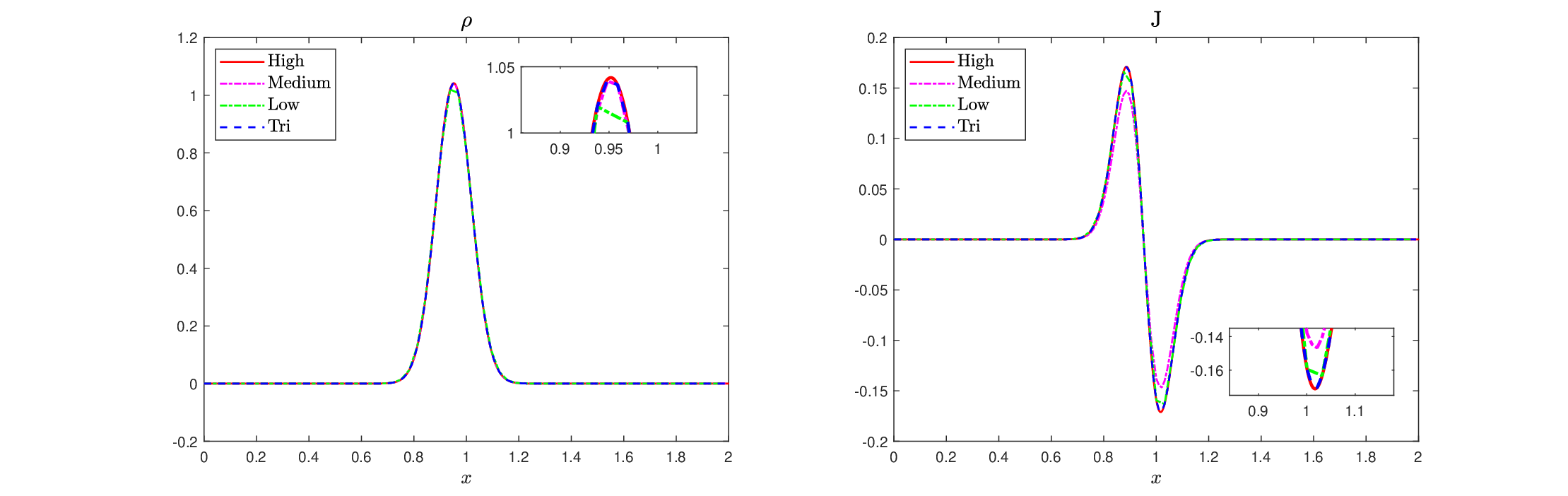}
	\includegraphics[width=0.49\textwidth]{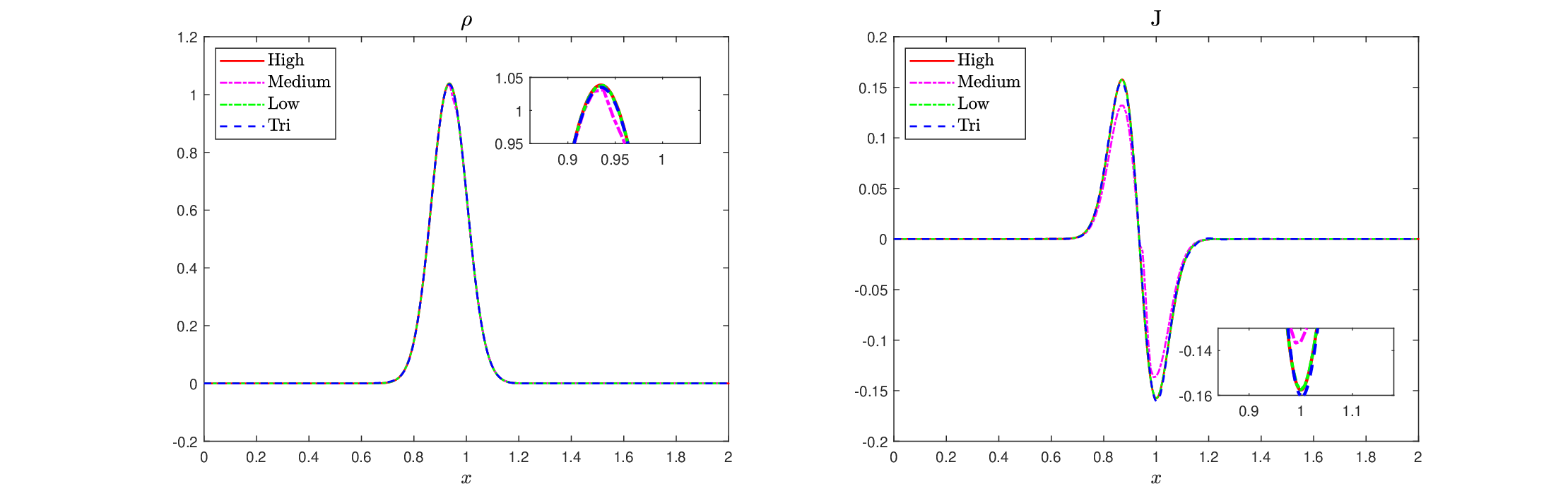}	\caption{Test 1. The mean and standard deviation of high- and tri-fidelity solutions $\rho$ and $J$ with $r=25$, at $t=0.01$, and comparison of the low-fidelity (FGA), medium-fidelity (LS), high-fidelity (TSFP) solution and the corresponding tri-fidelity approximations for fixed $z$ with different $\varepsilon$, i.e., $\varepsilon=1/64$ for the left two columns and $\varepsilon=1/256$ for the right two columns.} \label{Ex2Fig3N64}
\end{figure}

Regarding the tri-fidelity method, since the LS solver with finer meshes is more expensive than FGA, we choose FGA as the low-fidelity and LS as the medium-fidelity solver. The mean and standard deviation of the bi-fidelity solutions $\rho$ and $J$ are shown in Figure \ref{Ex2Fig3N64}  with $\varepsilon=1/64$ and $\varepsilon=1/256$ respectively, by evaluating the medium- and high- fidelity method only 25 times at $t=0.01$. Comparisons of the low-fidelity solution (FGA), medium-fidelity solution (LS), high-fidelity solution (TSFP) and the corresponding tri-fidelity approximations for an arbitrarily chosen sample point $z$ are also shown in Figure \ref{Ex2Fig3N64}. We can observe from the third row of Figure \ref{Ex2Fig3N64} that although the low- and medium-fidelity solutions are not accurate especially at positions where sharp transitions happen, yet the tri-fidelity approximation can still capture the solution well with random initial data.

\begin{figure}[!t]
	\centering
\includegraphics[width=0.38\textwidth]{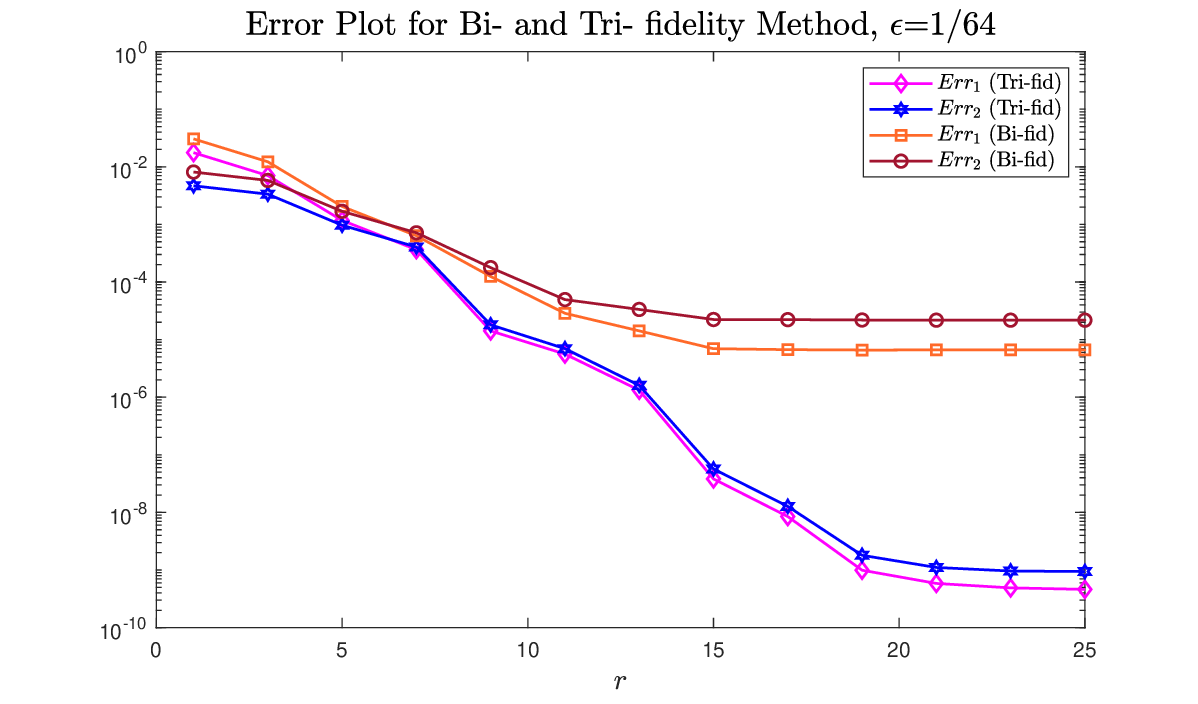}
\includegraphics[width=0.38\textwidth]{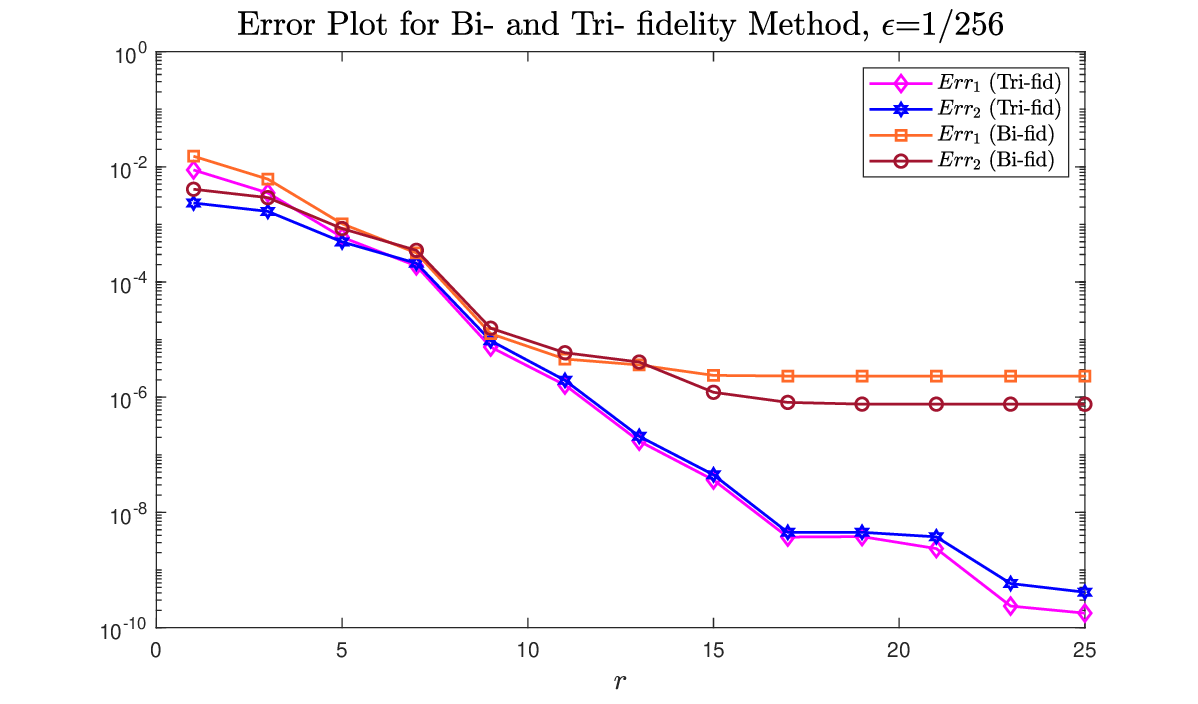}
	\caption{Test 1: Case I. The mean $L^2$ error of the bi-fidelity and tri-fidelity approximation $\rho, J$ with respect to the number of high-fidelity runs for different $\varepsilon$. }\label{Ex2Fig1Err2}
\end{figure}

Besides, we design another tri-fidelity method: based on the bi-fidelity method in Case I, keeping the low-fidelity (FGA) and high-fidelity (TSFP) solver unchanged, and using the TSFP method with the coarser grid ($N_{TSFP}^H = 10 N_{TSFP}^M$) as the medium-fidelity solver. We can see from Figure \ref{Ex2Fig1Err2} that the mean $L^2$ error of the tri-fidelity solution decays much more quickly than that of the bi-fidelity method for $r>10$. Figure \ref{Ex2Fig1Err2} shows that the tri-fidelity method performs better than the bi-fidelity method in this example.

\begin{figure}[!t]
	\centering
	\includegraphics[width=0.38\textwidth]{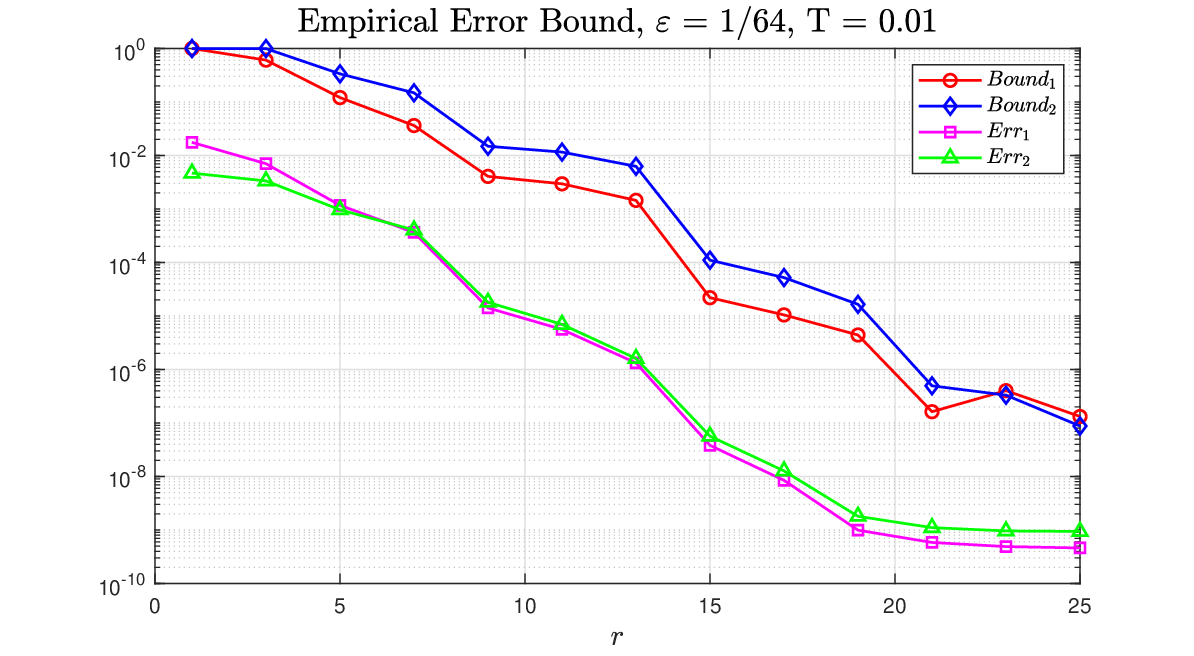}
	\includegraphics[width=0.38\textwidth]{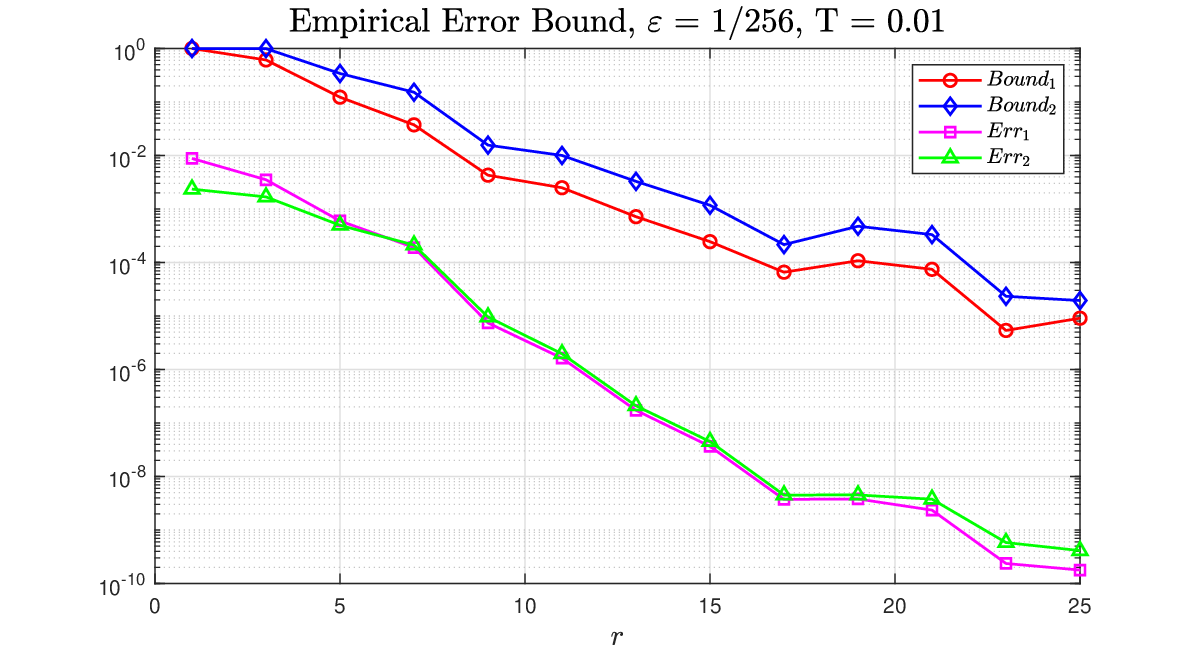}
	\caption{The mean $L^2$ error of the tri-fidelity approximation $\rho, J$, and the corresponding error bound estimation with respect to the number of high-fidelity runs for different $\varepsilon$.}\label{Fig_tri_bound}
\end{figure}

The mean $L^2$ error and corresponding error bound estimate with respect to the number of high-fidelity runs are shown in Figure \ref{Fig_tri_bound}. 
$Err_1$, $Err_2$ stand for the errors calculated by (4.1) for the macro quantities $\rho$ and $J$ respectively, while $Bound_1$, $Bound_2$ stand for the corresponding empirical error bounds computed by \eqref{empirical_bound} for the macro quantities $\rho$ and $J$ respectively.
It can be seen from Figure \ref{Fig_tri_bound} that the empirical error bounds (red and blue solid lines) are able to bound well the true $L^2$ errors (pink and green lines). By adding more high-fidelity solution samples one can reduce the tri-fidelity errors until they saturate, with the decay trend in consistent with the true $L^2$ errors and satisfy our expectations.


Below we make some remarks based on our experience of numerical simulations, between the medium- and low-fidelity solvers, we tend to choose the more expensive one as the medium-fidelity solver, which we only have to run $K$ ($K\ll N$) times. 
Since one needs to explore the low-fidelity solution space by repeatedly running the low-fidelity solvers $N$ times (with $N$ significantly large), the cheapest solver is always chosen as the low-fidelity model, regardless of its accuracy compared with medium-, not to say the high-fidelity solvers. 
Finally, it is noted that between the medium- and high-fidelity solvers, we always choose the more accurate one as the high-fidelity solver.

\subsection{Test 2}

In the first \textbf{Test 2(a)}, the initial condition is assumed by 
\begin{eqnarray}\label{test2}
	\begin{aligned}
		& n_0(x,z)=\exp\left(-\left(1+0.6\sum_{k=1}^{d_1} \frac{z_k^a}{2 k}\right)\left(x+1-0.4 \sum_{k=1}^{d_1} \frac{z_k^b}{2 k}\right)^2/\varepsilon\right), \quad\\
		& S_0(x)=x+1-0.4 \sum_{k=1}^{d_1} \frac{z_k^c}{2 k}, \quad x \in \mathbb{R} .	
	\end{aligned}
\end{eqnarray}
This example was based on \cite{BJM2002}, where we add uncertainties to the initial condition. Let $V(x)=x^2/2$, which is a harmonic oscillator.  
In Test 2(a), the spatial domain is chosen to be $[-\pi,\pi]$, and periodic boundary condition is assumed. 
Here $\mathbf{z}^a=\left(z_1^a, \cdots, z_{d_1}^a\right), \mathbf{z}^b=\left(z_1^b, \cdots, z_{d_1}^b\right)$ and $\mathbf{z}^c=\left(z_1^c, \cdots, z_{d_1}^c\right)$ represent the random variables in initial condition of $n_0(x)$ and $S_0(x)$. Set $d_1=5$, thus this is a $d=20$ dimensional problem in the random space. 

We first choose FGA as the low-fidelity solver, TSFP as the high-fidelity solver and investigate the long-time behavior of the bi-fidelity solution.
We consider the two different regimes with $\varepsilon$ varying from $\varepsilon=1/32$ to $\varepsilon=1/128$ in this test. When we choose $\varepsilon=1/32$, we take $N_{FGA} = 384$ and $N_{TSFP} =10 N_{FGA}$ grid points for FGA and TSFP method. 
When $\varepsilon=1/128$, for spatial discretization, we take $N_{FGA} = 1536$ and $N_{TSFP} =5 N_{FGA}$ grid points for FGA and TSFP method. 
Let $\Delta t=\Delta x/20$ in both high- and low-fidelity models, with the final time $T=6.0$.
We choose 800 number of particles when $\varepsilon=1/32$ and 1300 particles when $\varepsilon=1/128$ for the FGA method.

\begin{figure}[!t]
\centering
\includegraphics[width=0.38\textwidth]{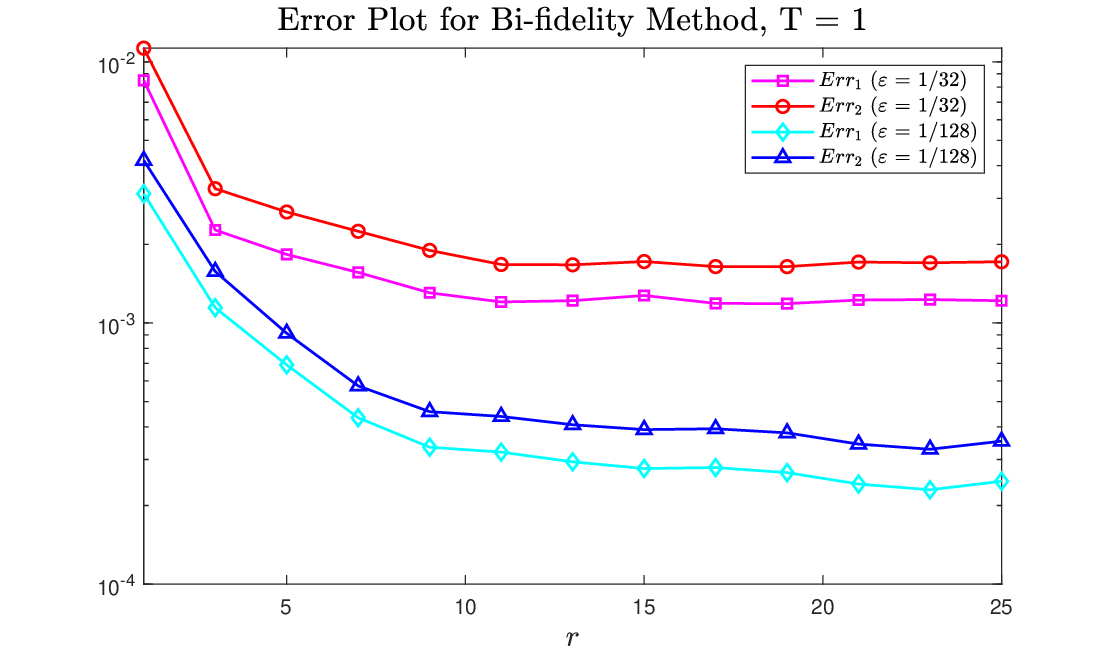}
\includegraphics[width=0.38\textwidth]{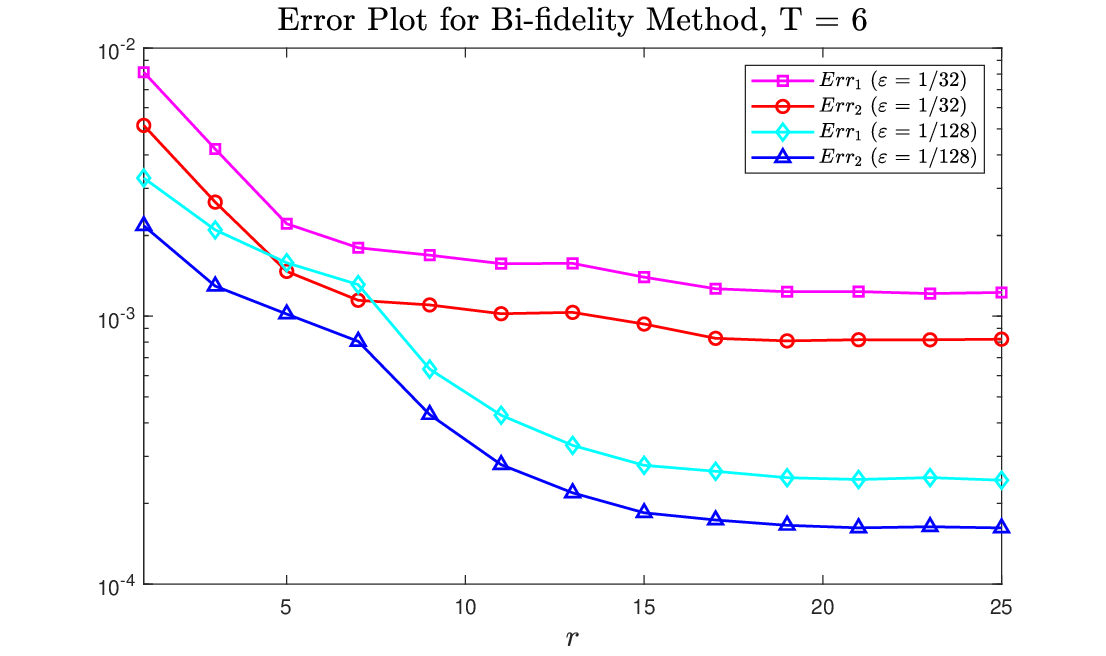}
\caption{Test 2. The mean $L^2$ errors of the bi-fidelity approximation $\rho, J$ with respect to the number of high-fidelity runs for different $\varepsilon$. Here, $Err_1, Err_2$ stand for the errors of $\rho$ and $J$.  }\label{Ex1Fig1Err}
\end{figure}

\begin{figure}[!t]
	\centering
	\includegraphics[width=0.38\textwidth]{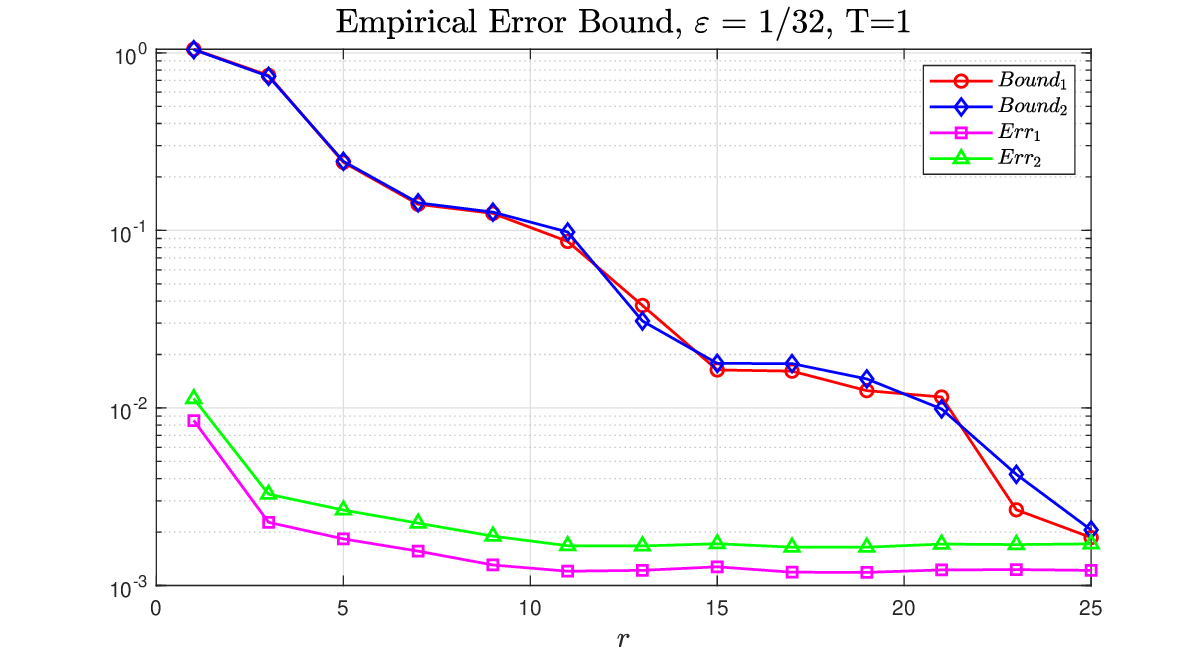}
	\includegraphics[width=0.38\textwidth]{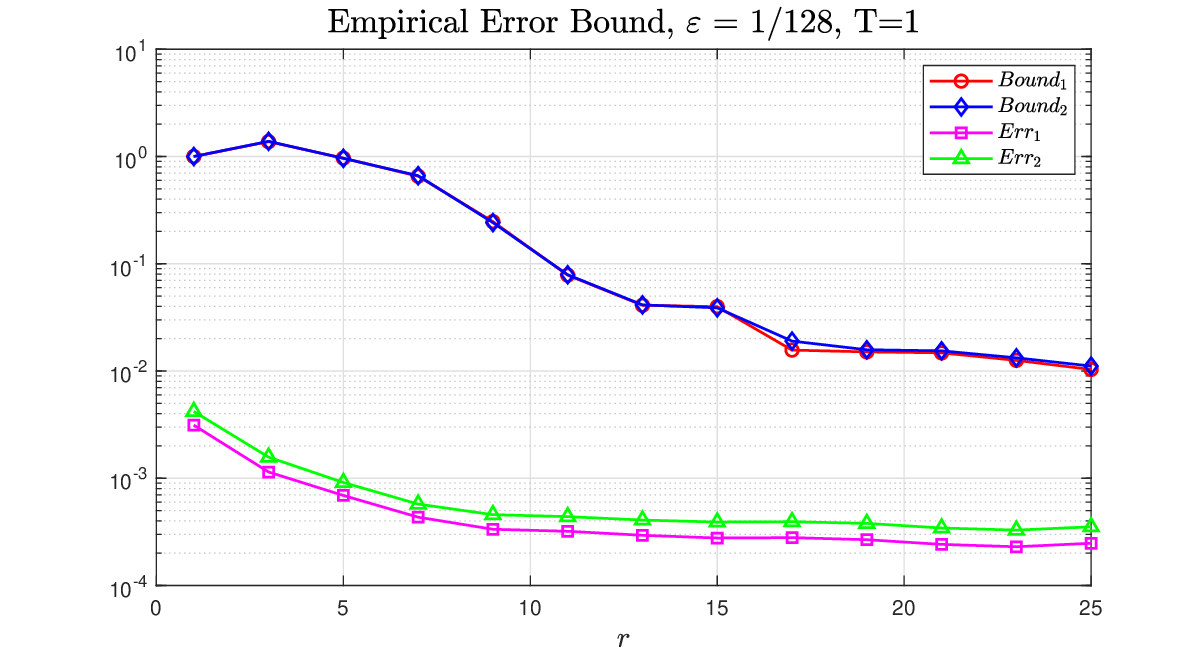}
	\caption{The mean $L^2$ error of the bi-fidelity approximation $\rho, J$, and the corresponding error bound estimation with respect to the number of high-fidelity runs for different $\varepsilon$. }\label{Fig_bi_bound}
\end{figure}

Figure \ref{Ex1Fig1Err} shows the mean $L^2$ errors of $\rho$ and $J$ between the high- and bi-fidelity solutions with different $\varepsilon$. It is clear that the error decays fast with the number of high-fidelity runs. In addition, when $\varepsilon$ decreases, the error between the high- and bi-fidelity solutions grows smaller. 
From Figure \ref{Fig_bi_bound}, one can see that the practical error bound estimators bound well the true bi-fidelity errors.

\begin{figure}[!t]
	\centering
	\includegraphics[width=0.49\textwidth]{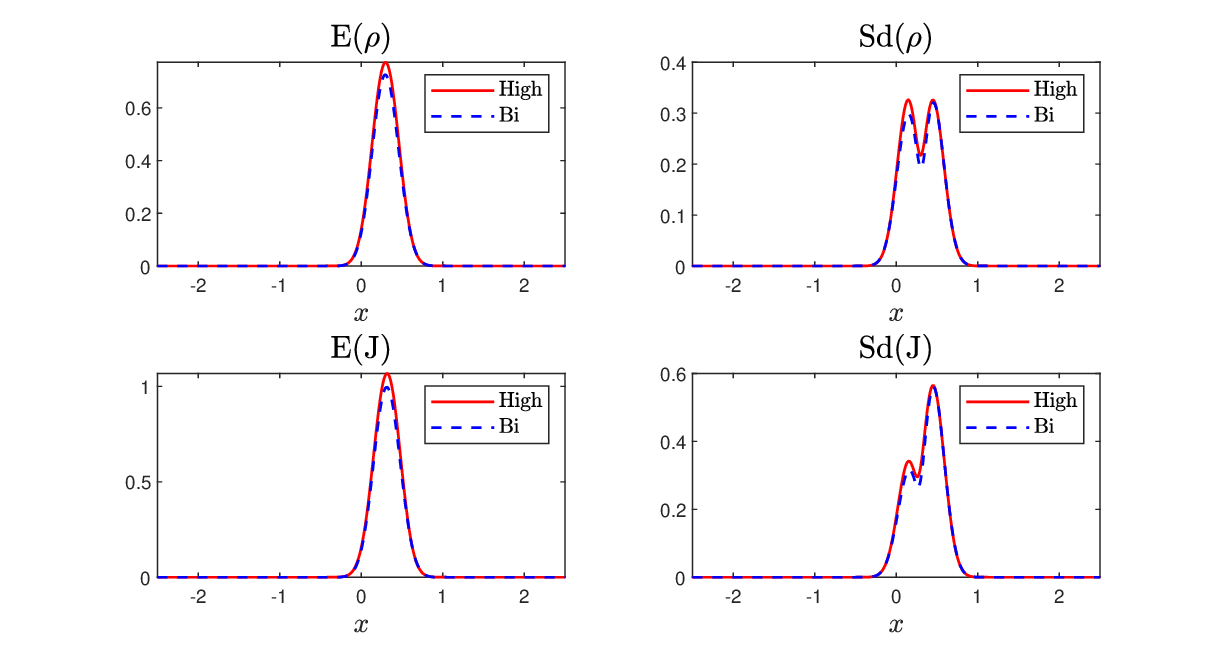}
	\includegraphics[width=0.49\textwidth]{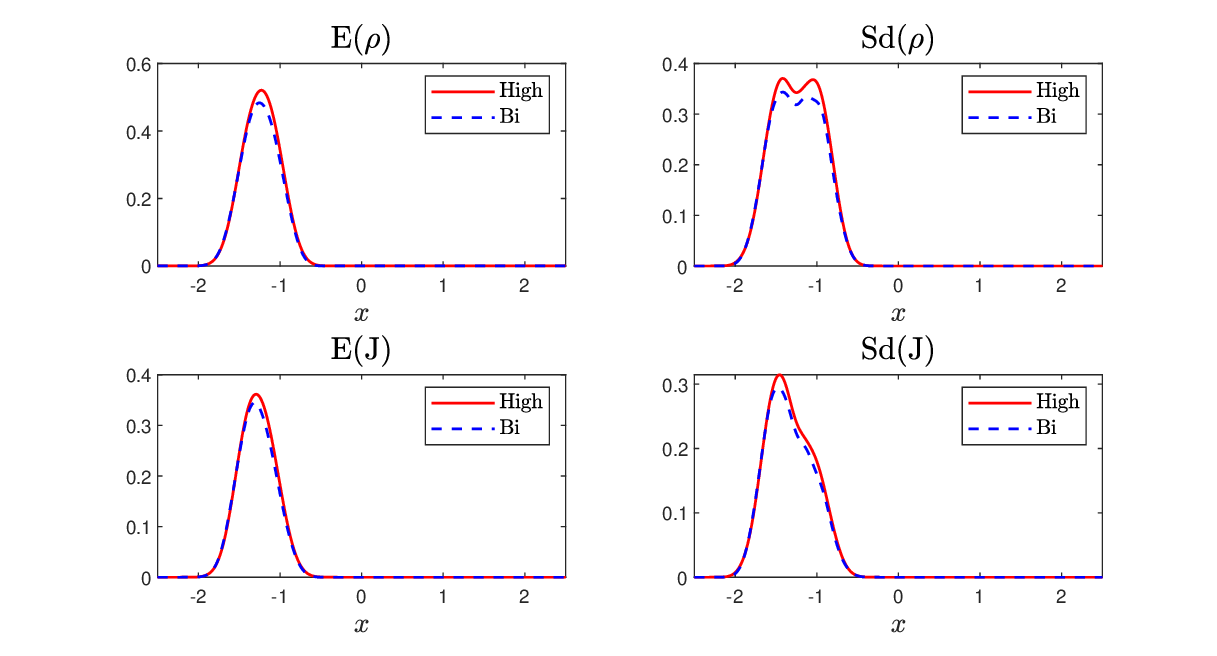}
	\caption{Test 2. The mean and standard deviation of high- and bi-fidelity solutions $\rho$ and $J$ with $r=25, \varepsilon=1/32$ at different time, with $t=1$ for the left two columns, 
		and $t=6$ for the right two columns.}\label{Ex1Fig2N32}
\end{figure}

\begin{figure}[!t]
	\centering
	\includegraphics[width=0.49\textwidth]{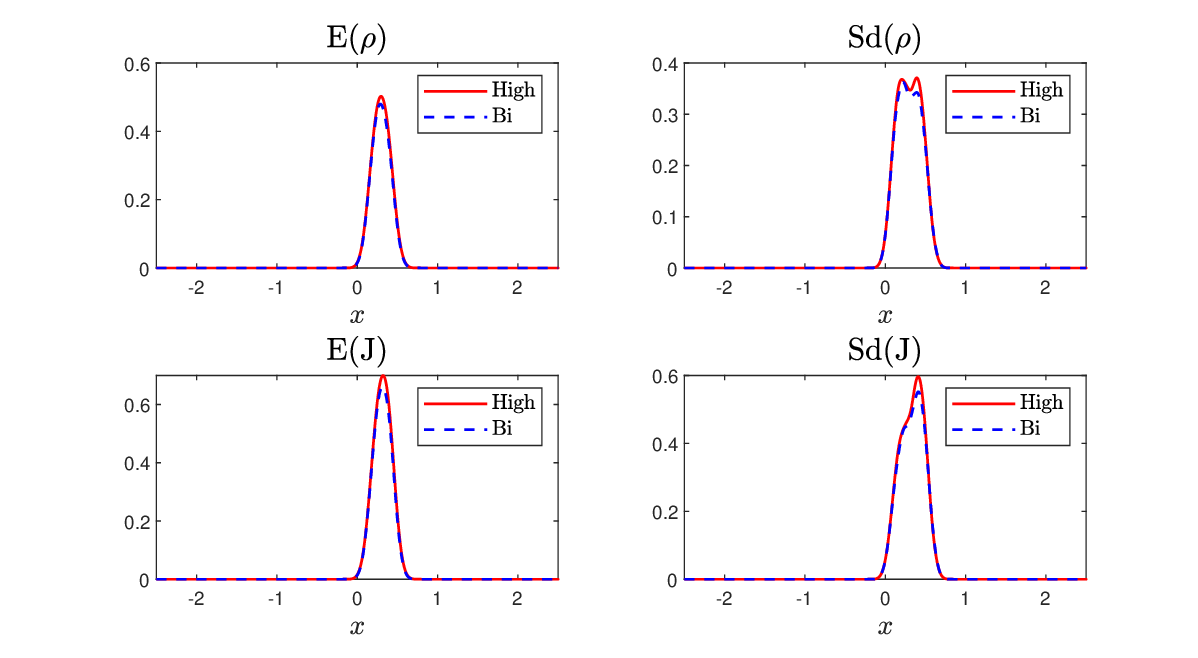}
	\includegraphics[width=0.49\textwidth]{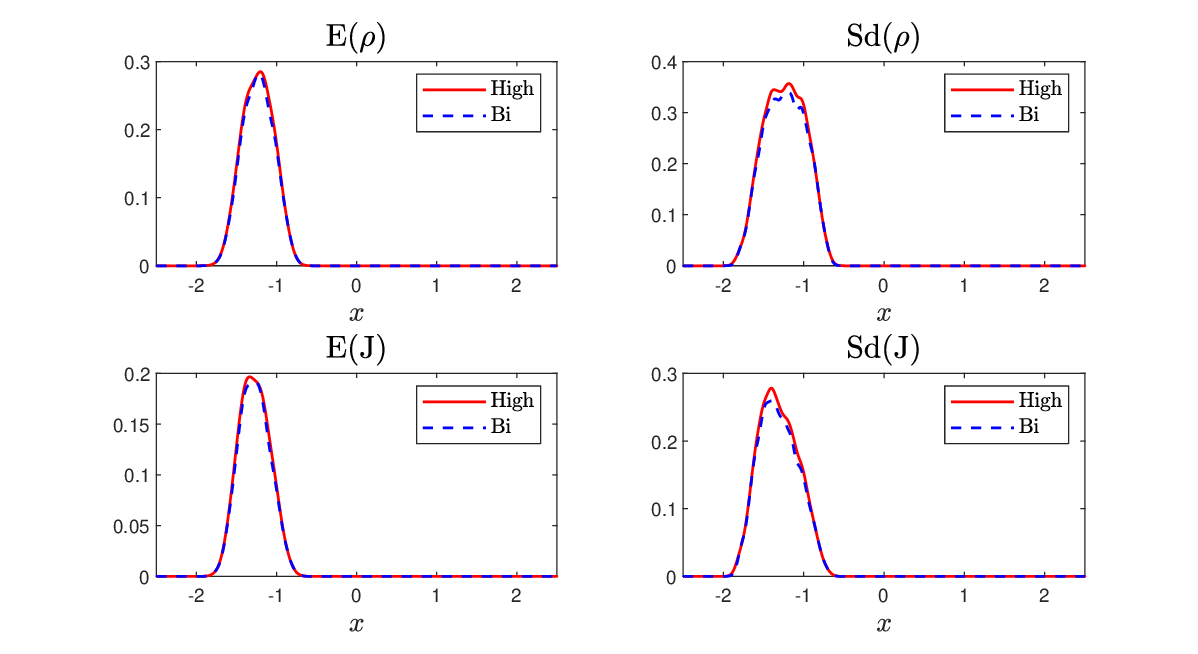}
	\caption{Test 2. The mean and standard deviation of high- and bi- fidelity solutions $\rho$ and $J$ with $r=25, \varepsilon=1/128$ at different time, with $t=1$ for the left two columns, 
		and $t=6$ for the right two columns.}\label{Ex1Fig2N128}
\end{figure}

In this test, the mean and standard deviation of the bi-fidelity solutions $\rho$ and $J$ are shown in Figure \ref{Ex1Fig2N32} with $\varepsilon=1/32$ and Figure \ref{Ex1Fig2N128} with $\varepsilon=1/128$ by adopting the high-fidelity solver only 25 times at both $t=1$ and $t=6$. We can see that all the mean and standard deviation of the bi-fidelity approximation of $\rho$ and $J$ match well with the high-fidelity solutions by using 25 high-fidelity runs. When $\varepsilon=1/32$, the high-fidelity model (TSFP) costs approximately 7 times of the low-fidelity solver (FGA) with the final time $T=6$ (the former takes 15.64 seconds while the latter takes 2.03 seconds for one single run). Therefore, a significant speedup is quite noticeable in this case.
This shows the efficiency of our method to capture the long time behavior of the solution.

\begin{figure}[!t]
	\centering
	\includegraphics[width=0.49\textwidth]{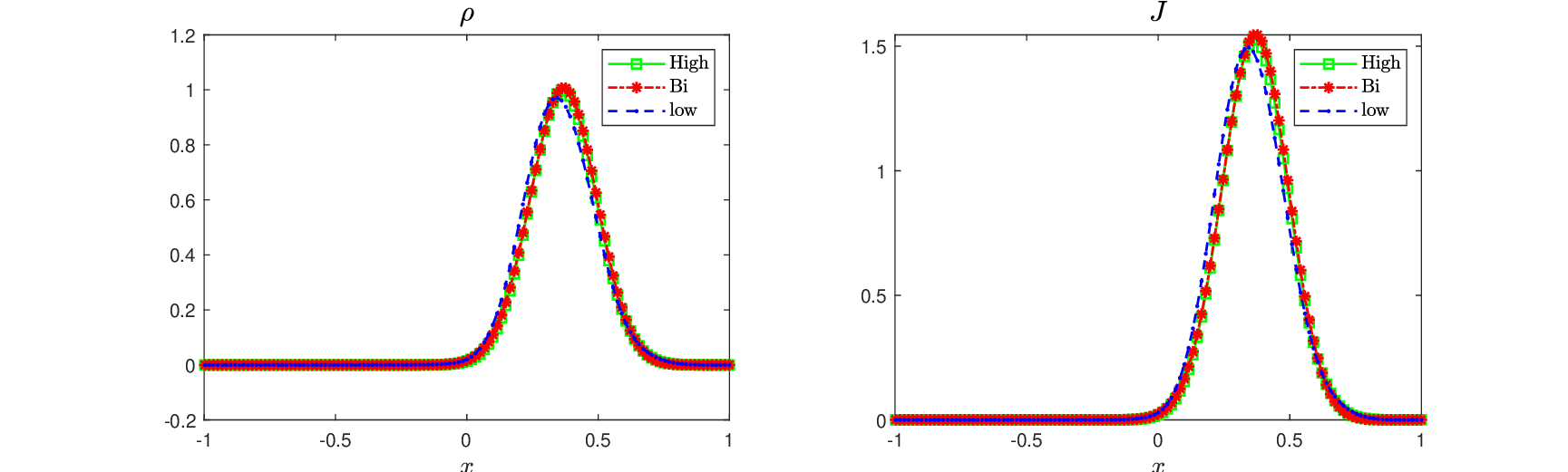}
	\includegraphics[width=0.49\textwidth]{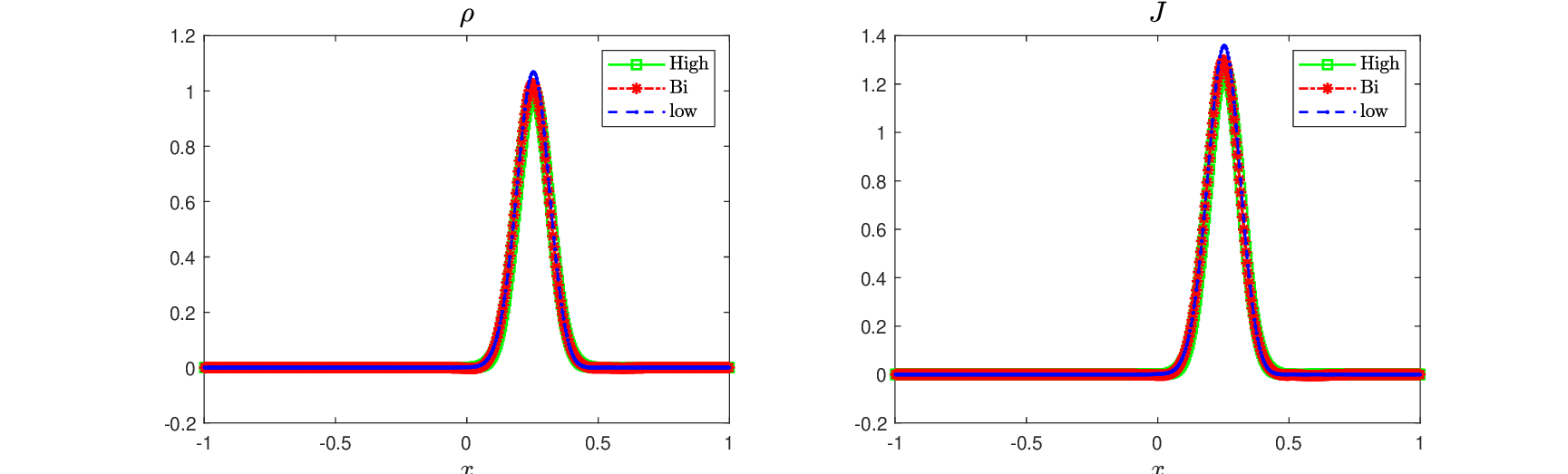}
	\caption{Test 2.  Comparison of the low-fidelity solution (FGA), high-fidelity solution (TSFP) and the corresponding bi-fidelity approximations with $r=25$, $t=1.0$ for a fixed $z$, $\varepsilon=1/32$ for the left two figures and $\varepsilon=1/128$ for the right two figures.   }\label{Ex1Fig3N128}
\end{figure}
	
To further illustrate the performance of our bi-fidelity method, we compare the high-fidelity, low- and the corresponding bi-fidelity solutions (with $r=25$) for a particular sample point $z$. One observes from  Figure \ref{Ex1Fig3N128} ($\varepsilon=1/128$) that the high- and bi-fidelity solutions match really well, whereas the low-fidelity solutions appear inaccurate at some points in the spatial domain.

\begin{figure}[!t]
	\centering
	\includegraphics[width=0.49\textwidth]{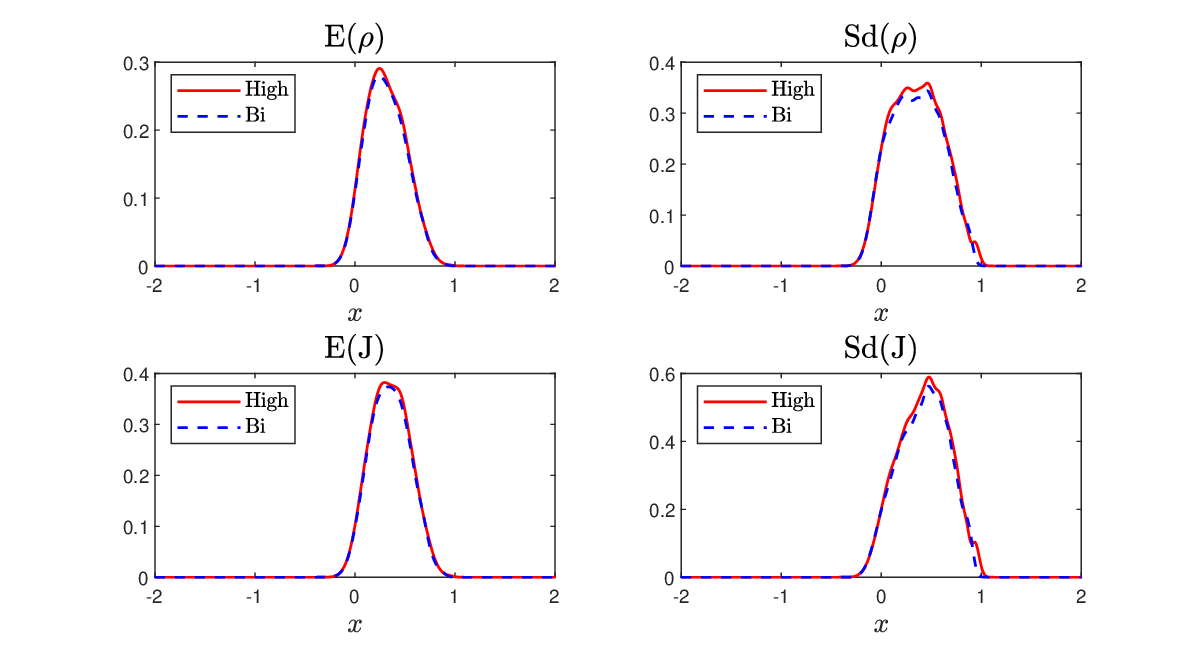}
	\caption{Test 2. The mean and standard deviation of high- and bi- fidelity solutions $\rho$ and $J$ with $r=25, \varepsilon=1/128$ at $t=1$. }\label{Ex1Fig2BothRand}
\end{figure}

\begin{figure}[!t]
	\centering\label{Randxzfig}
	\includegraphics[width=0.49\textwidth]{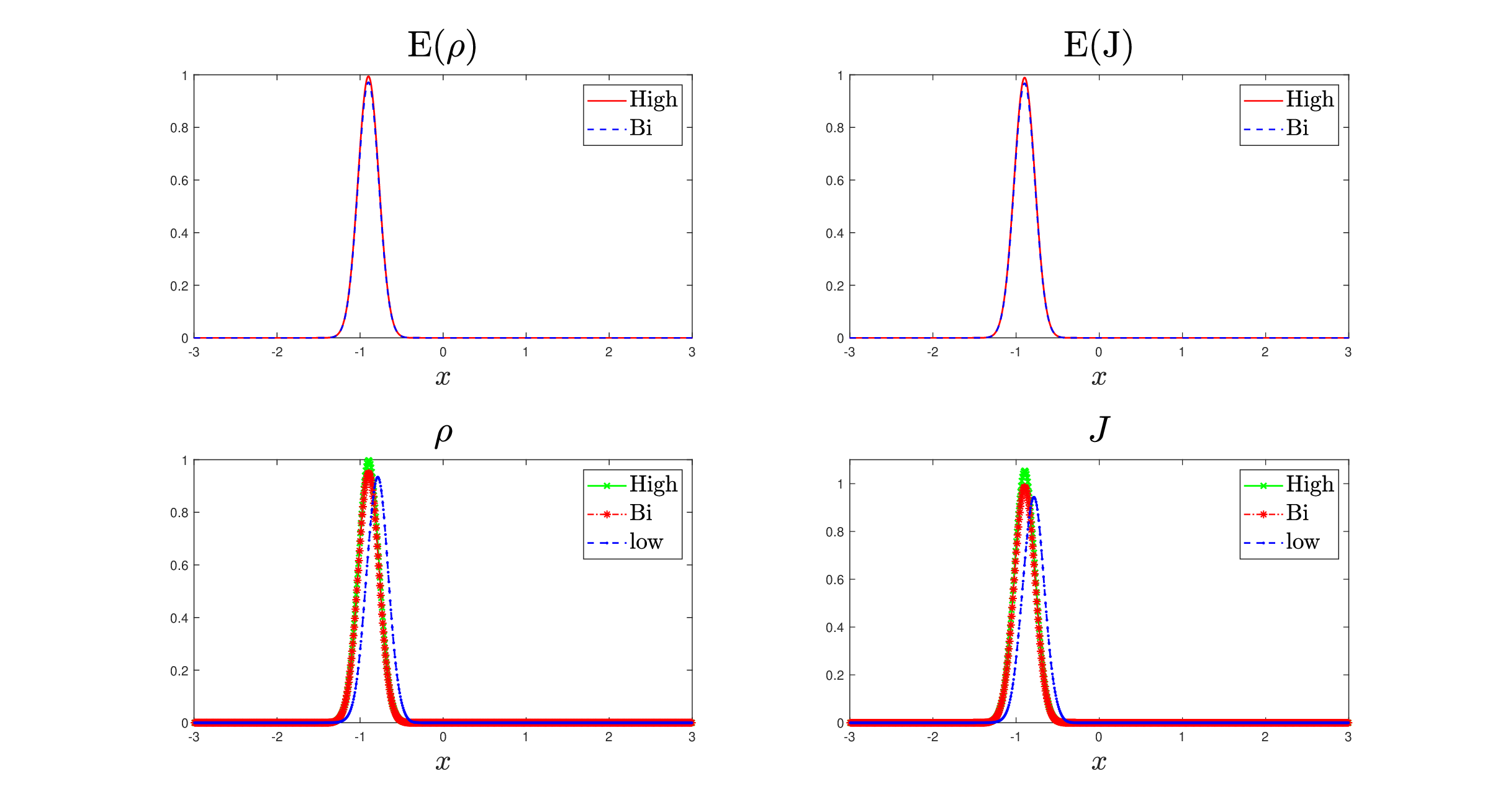}
	\caption{Expectation of high-and bi-fidelity solutions $\rho$ and $J$ with $r=25, \varepsilon=1 / 32$ (fist row) and comparison of the low-fidelity solution (FGA), high-fidelity solution (TSFP) and the corresponding bi-fidelity approximations for a fixed $z$ (second row).}
\end{figure}

\textbf{Test 2(b)}. In this test, we consider the random uncertainties coming from both the potential function $V$ and the initial data. 
Figure \ref{Ex1Fig2BothRand} shows the mean and standard deviation of high- and bi-fidelity solutions $\rho$ and $J$ at $r=25, \varepsilon=1/128, t=1$ with a simple random potential $V(x,z) = x^2/2 + 0.5 \sum_{k=1}^{d_1} \frac{z_k}{2k}$ together with the random initial condition \eqref{test2}. Considering  a more complex random potential $V(x,z) = \sum_{k=1}^{d_1} z_k\cdot (kx/10)^2 $, Figure \ref{Randxzfig} ($\varepsilon=1/32$) shows the expectation of high- and bi-fidelity solutions and comparison of the low-fidelity solution, high-fidelity solution and the corresponding bi-fidelity approximations for a particular sample point. The efficiency of our bi-fidelity method has been validated.

\begin{figure}[!t]
	\centering
	\includegraphics[width=0.47\textwidth]{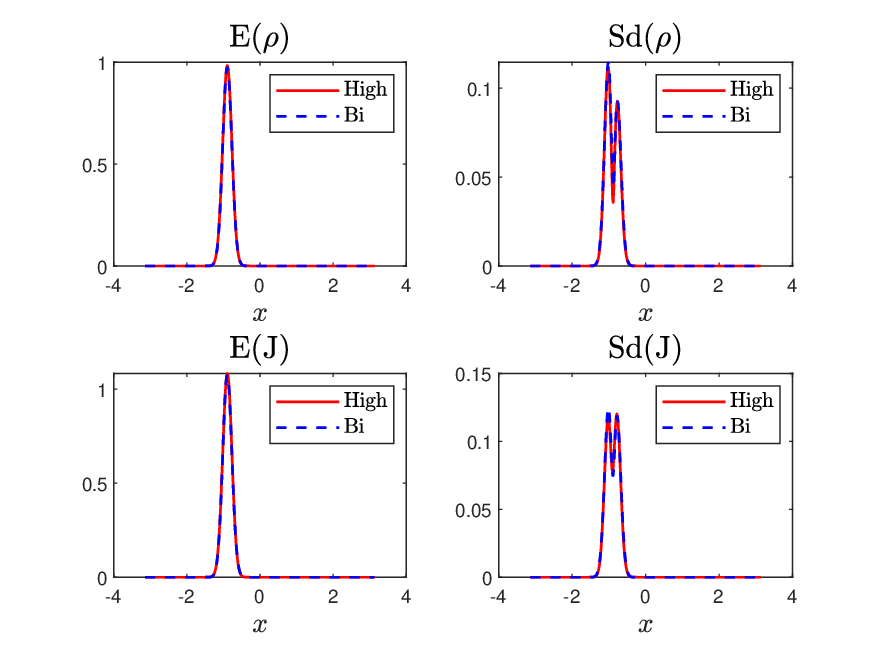}
	\includegraphics[width=0.47\textwidth]{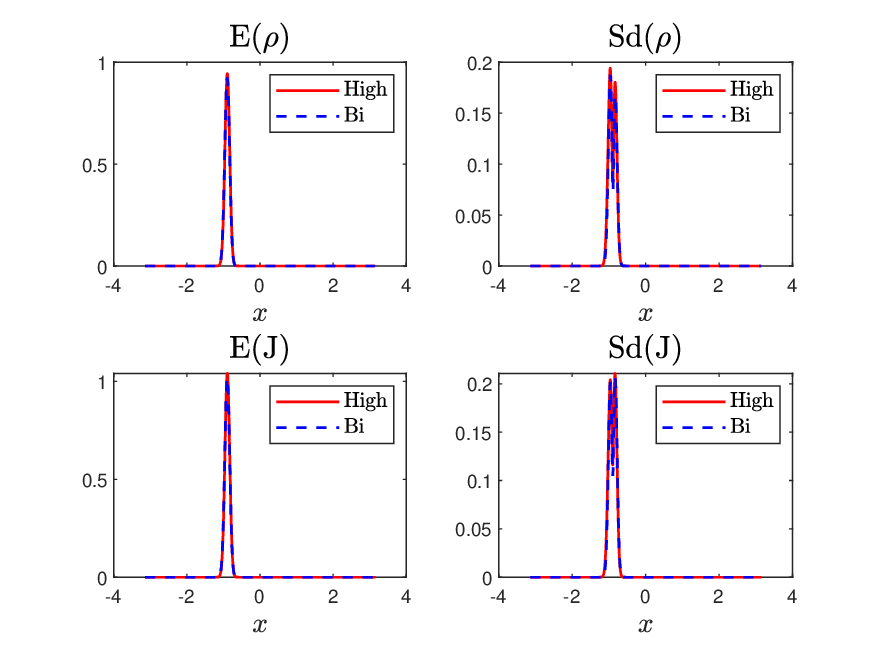}
    	\includegraphics[width=0.47\textwidth]{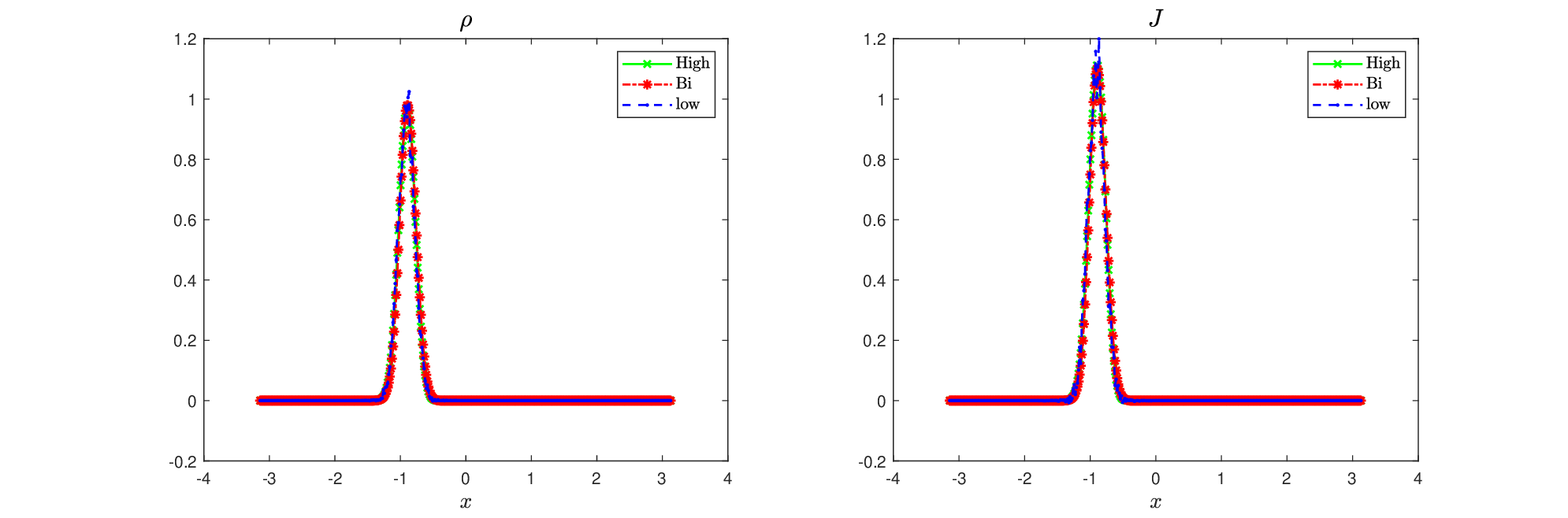}
        \includegraphics[width=0.48\textwidth]{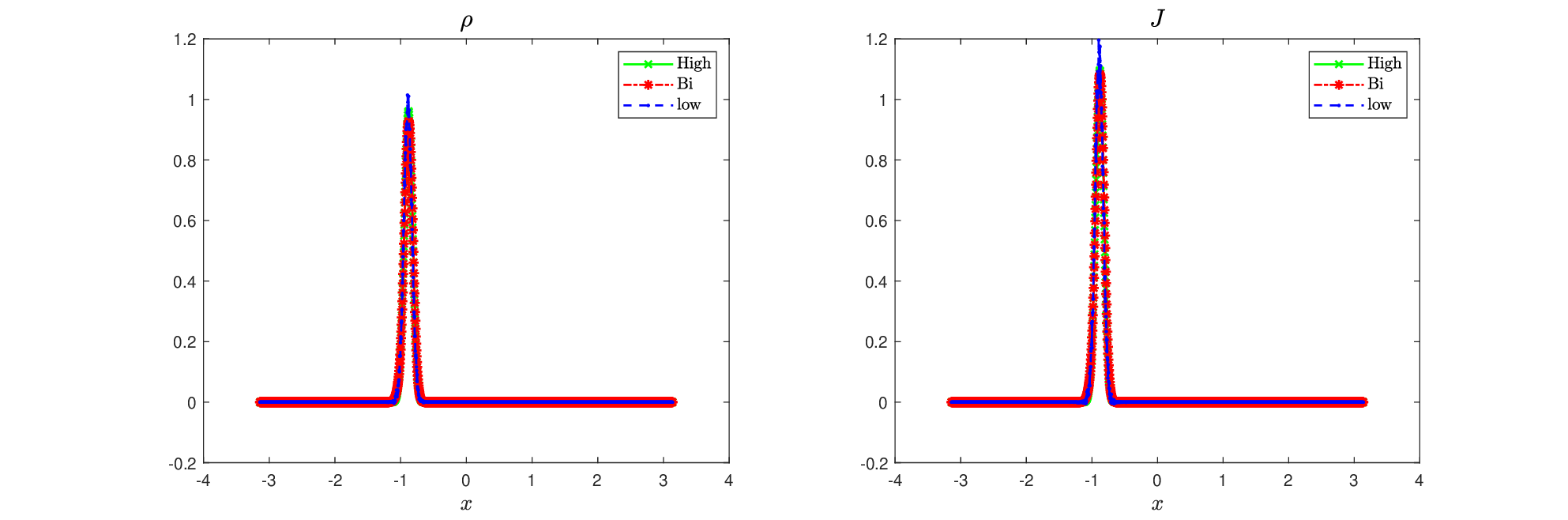}
	\caption{Upper two lines are the mean and standard deviation of high- and bi- fidelity solutions $\rho$
and $J$ with $r = 25$ for different $\varepsilon$, that is, $\varepsilon=1/32$ for left four figures and $\varepsilon=1/128$ for right four figures. Last line is comparison of the low-fidelity solution (FGA), high-fidelity solution (TSFP)
and the corresponding bi-fidelity approximations with $r = 25$, $\varepsilon=1/32$ (left two figures) and  $\varepsilon=1/128$ (right two figures) for a fixed $z$.
}\label{Fig_KL}
\end{figure}

\textbf{Test 2(c).}
In this final test, uncertain parameters depending on the spatial variable are studied in order to better illustrate the performance of our method. We believe these random parameters are more complex and closer to real applications. Here we use the widely used Karhunen-Lo\`eve (KL) expansion \cite{KL}
to represent the randomness of the initial condition.
 Specifically, the following initial condition by KL expansion is assumed by
\begin{align*}
    {\small \begin{aligned}
& n_0(x, z)=\exp \left(-\left(1+0.6 \sum_{k=1}^{d_1} \sqrt{\lambda_k} z_k^a\varphi_k(x)\right)\left(x+1-0.4 \sum_{k=1}^{d_1} \sqrt{\lambda_k} z_k^b\varphi_k(x)\right)^2 / \varepsilon\right), \\
& S_0(x)=x+1-0.4 \sum_{k=1}^{d_1} \sqrt{\lambda_k} z_k^c\varphi_k(x), \quad x \in \mathbb{R},
\end{aligned}}
\end{align*}
where $z_k^a, z_k^b, z_k^c$ are i.i.d. random variables with zero mean and unit covariance, $\{\varphi_k,\lambda_k\}_{k=0}^\infty$ are orthonormal eigenfunction and corresponding eigenvalues of the covariance operator $K_f$
$$K \phi(x):=\int_{-\pi}^{\pi} c(x-y) \phi(y) \mathrm{d} y, \quad c(x-y)=\sigma^2 e^{-|x-y|^2 /l^2},$$
where $\sigma$ is the root mean square height of the profile  and $l$ is the correlation length.
Since the eigenvalues are in a descending order, one only need the first several terms to
represent the random initial input.
In this test, we choose $l=0.5$ and $\sigma=0.05$. The mean and standard deviation of the bi-fidelity solutions $\rho$ and $J$ are shown in the upper two lines in Figure \ref{Fig_KL} with $\varepsilon=1 / 32, 1 / 128$ by adopting the high-fidelity solver only 25 times at $t=0.1$. We can see that all the mean and standard deviation of the bi-fidelity approximation of $\rho$ and $J$ match well with the high-fidelity solutions by using 25 high-fidelity runs. One observes from the last line in Figure \ref{Fig_KL} that the high- and bi-fidelity solutions match really well, whereas the low-fidelity solutions are not accurate at some points in the spatial domain.

\section{Conclusion}
\label{sec:conclusion}

In this work, we study robust bi-fidelity and tri-fidelity method for the Schr\"odinger equation with random parameters. 
Different fidelity methods for solving the Schr\"odinger equation including the time-splitting Fourier pseudo-spectral (TSFP) method, the frozen Gaussian Approximation (FGA), and level set (LS) method for the semi-classical limit of the Schr\"odinger equation are explored. Based on the error estimates for different numerical schemes, we conduct an error analysis of the bi-fidelity method together with implementation on the empirical error bound estimations. 
Moreover, we compare the bi-fidelity and tri-fidelity approximations, which will hopefully provide us some insight on appropriate choices of low- and medium-fidelity models for more general UQ problems. Extensive numerical examples have demonstrated the accuracy and efficiency of our proposed bi-fidelity and tri-fidelity method for solving the Schr\"odinger equation with random inputs. In future work, we plan to investigate the error analysis for multi-fidelity methods and develop efficient multi-fidelity solvers for the Schr\"odinger equation with more complex and higher-dimensional random parameters.

\section*{Acknowledgement}
L. Liu acknowledges the support by National Key R\&D Program of China (2021YFA1001200), Ministry of Science and Technology in China, Early Career Scheme (24301021) and General Research Fund (14303022 \& 14301423) funded by Research Grants Council of Hong Kong from 2021-2023.  
Y. Lin was partially supported by National Key R\&D Program of China (2024YFA1016100) and National Natural Science Foundation of China (12201404).
The computations in this paper were run on the Siyuan-1 cluster supported by the Center for High Performance Computing at Shanghai Jiao Tong University.
We thank Prof. Shi Jin and reviewers for their helpful suggestions on improving the manuscript. 

\appendix
\section{Proof of Lemma \ref{thm:L2TSFP}} \label{app:A}
	Define 
	$$
	Y_N:=\operatorname{span}\left\{e^{i \mu_l(x-a)}, x \in \bar{\Omega}, l \in \mathcal{T}_N\right\}, \mathcal{T}_N=\left\{l \mid l=-\frac{N}{2}, \ldots, \frac{N}{2}-1\right\},
	$$
	where $\mu_l=\frac{2 \pi l}{b-a}$. Let $P_N$ be the standard $L^2$-projection operator onto $Y_N$, and  $I_N$ be the trigonometric interpolation operator \cite{Shen2011spectral}, i.e.,
	$$
	P_N u=\sum_{l \in \mathcal{T}_N} \widehat{u}_l e^{i \mu_l(x-a)}, \quad I_N u=\sum_{l \in \mathcal{T}_N} \widetilde{u}_l e^{i \mu_l(x-a)}, \quad x \in \bar{\Omega},
	$$
	where
	$
	\widehat{u}_l=\frac{1}{b-a} \int_a^b u(x) e^{-i \mu_l(x-a)} d x, \quad \widetilde{u}_l=\frac{1}{N} \sum_{j=0}^{N-1} u(x_j) e^{-i \mu_l\left(x_j-a\right)}, \quad l \in \mathcal{T}_N.
	$
	Since
	$$
	I_N \psi^n-\psi\left(t_n\right)=I_N \psi^n-P_N\left(\psi\left(t_n\right)\right)+P_N\left(\psi\left(t_n\right)\right)-\psi\left(t_n\right),
	$$
	under Assumption \ref{AssumpAB}, the standard Fourier projection properties \cite{Shen2011spectral,BCF2023}  yield
	\begin{align}\label{projectionErr}
		\left\|I_N \psi^n-\psi\left(t_n\right)\right\|_{L^2} \leq\left\|I_N \psi^n-P_N\left(\psi\left(t_n\right)\right)\right\|_{L^2}+C_1 h^m, .
	\end{align}
	Consider the error function $e^n$ at $t_n$ as
	$
	e^n:=e^n(x)=I_N \psi^n-P_N \psi\left(t_n\right) ~ (n \geq 1),
	$
	and $\left\|e^0\right\|_{L^2} \leq C_2 h^m$ implied by the standard projection and interpolation results. We can write the error function $e^n$ as
	\begin{align}\label{eDef}
		e^{n+1}=I_N \psi^{n+1}-P_N \mathcal{S}_\tau\left(P_N \psi\left(t_n\right)\right)+\mathcal{E}^n,
	\end{align}
	with the local truncation error
	$\mathcal{E}^n :=\mathcal{E}^n(x)=P_N \mathcal{S}_\tau\left(P_N \psi\left(t_n\right)\right)-P_N \psi\left(t_{n+1}\right)$. Here $\mathcal{S}_\tau$ is defined by
 \begin{align*}\label{def:S_tau}
     \mathcal{S}_\tau\left(P_N \psi\left(t_n\right)\right)=e^{i \frac{\e}{2}\frac{\tau}{2} \Delta} e^{-i  \frac{\tau}{\varepsilon} V(x)} e^{i \frac{\e}{2}\frac{\tau}{2} \Delta} P_N \psi\left(t_n\right), \quad x\in\bar{\Omega}.
 \end{align*}
	Taking the $L^2$-norm on both sides of \eqref{eDef}, one has
	$$	\begin{aligned}
		\left\|e^{n+1}\right\|_{L^2} & \leq\left\|I_N \psi^{n+1}-P_N \mathcal{S}_\tau\left(P_N \psi\left(t_n\right)\right)\right\|_{L^2}+\left\|\mathcal{E}^n\right\|_{L^2}.
	\end{aligned} $$
	
	\textbf{Step 1. Estimate for $||I_N \psi^{n+1}-P_N \mathcal{S}_\tau\left(P_N \psi\left(t_n\right)\right)||$.}
	The fully discrete scheme \eqref{TSFP}  can be written as
	{\small $$
	\begin{aligned}
			& I_N \psi^{n+1}=e^{i \frac{\varepsilon}{2} \frac{\tau}{2} \Delta}\left(I_N \psi^{(**)}\right), I_N\left(\psi^{(**)}\right)=I_N\left(e^{-i \frac{\tau}{\varepsilon} V(x)} \psi^{(*)}\right), I_N \psi^{(*)}=e^{i \frac{\varepsilon}{2} \frac{\tau}{2} \Delta} I_N \psi^n, \\
			& P_N\left(\mathcal{S}_\tau\left(\psi\left(t_n\right)\right)\right)=e^{i \frac{\varepsilon}{2}\frac{\tau}{2} \Delta}\left(P_N \psi^{\langle **\rangle}\right), \psi^{\langle **\rangle}=e^{-i \frac{\tau}{\varepsilon} V(x)} \psi^{\langle *\rangle}, \psi^{\langle *\rangle}=e^{i\frac{\varepsilon}{2} \frac{\tau}{2} \Delta} P_N \psi\left(t_n\right).
	\end{aligned}  
	$$}
	Since $I_N$ and $P_N$ are identical on $Y_N$ and $e^{i \frac{\varepsilon}{2} \frac{\tau}{2} \Delta }$ preserves the $H^k$-norm $(k \geq 0)$, using Taylor expansion $e^{-i \frac{\tau}{\varepsilon} V(x)}=1-i \frac{\tau}{\varepsilon} V(x) \int_0^1 e^{-i  \theta \frac{\tau}{\varepsilon} V(x)} d \theta$ and Assumption \ref{AssumpAB}, one has
	\begin{align}\label{2.20}
		\begin{aligned}
			&\left\|I_N \psi^{n+1}-P_N \mathcal{S}_\tau(P_N \psi(t_n))\right\|_{L^2}			=\left\|I_N \psi^{(**)}-P_N \psi^{\langle **\rangle}\right\|_{L^2} \\
   \leq& \left\|I_N \psi^{(**)}-I_N \psi^{\langle **\rangle}\right\|_{L^2} + \left\|I_N \psi^{\langle **\rangle}-P_N \psi^{\langle **\rangle}\right\|_{L^2},
		\end{aligned} 
	\end{align} 
 and
	\begin{align}\label{2.21}
		\begin{aligned}
			\|P_N \psi^{\langle **\rangle}-I_N \psi^{\langle **\rangle}\|_{L^2}=&\left\|\frac{\tau}{\varepsilon}(P_N-I_N)\left(V(x) \int_0^1 e^{-i  \theta \frac{\tau}{\varepsilon} V(x)} d \theta \psi^{\langle *\rangle}\right)\right\|_{L^2} \\
			\leq& C_3 \frac{\tau}{\varepsilon} h^m,
		\end{aligned}
	\end{align} 
	where $C_3$ is obtained from Fourier interpolation and projection properties together with $\left\|V(x) \int_0^1 e^{-i  \theta \frac{\tau}{\varepsilon} V(x)} d \theta \psi^{\langle *\rangle}\right\|_{L^2} \lesssim \|V\|_{L^2}\left\|\psi\left(t_n\right)\right\|_{L^2}$. In addition, by direct computation and Parseval's identity, we can derive
	\begin{align}\label{2.22}
		\begin{aligned}
			\|I_N \psi^{(**)}-I_N \psi^{\langle **\rangle}\|_{L^2} & =\sqrt{h \sum_{j=0}^{N-1}\left|\psi_j^{(**)}-\psi^{\langle **\rangle}\left(x_j\right)\right|^2}=\sqrt{h \sum_{j=0}^{N-1}\left|\psi_j^{(1)}-\psi^{\langle *\rangle}\left(x_j\right)\right|^2} \\
			& =\left\|I_N \psi^{(*)}-I_N \psi^{\langle *\rangle}\right\|_{L^2}=\left\|I_N \psi^n-P_N \psi\left(t_n\right)\right\|_{L^2} 
			 =\left\|e^n\right\|_{L^2} .
		\end{aligned}
	\end{align}
	Combining  \eqref{2.20}, \eqref{2.21} and  \eqref{2.22} together, we obtain 
	\begin{align}\label{step1}
		\begin{aligned}
			&\left\|I_N \psi^{n+1}-P_N \mathcal{S}_\tau\left(P_N \psi\left(t_n\right)\right)\right\|_{L^2}\\
			 \leq &\left\|I_N \psi^{(**)}-I_N \psi^{\langle **\rangle}\right\|_{L^2} + \left\|I_N \psi^{\langle **\rangle}-P_N \psi^{\langle **\rangle}\right\|_{L^2}\\
			 \leq &\left\|e^n\right\|_{L^2} +  C_3 \frac{\tau}{\varepsilon} h^m,
		\end{aligned}
	\end{align} 
	
	\textbf{Step 2. Estimate for $||\mathcal{E}^n||$.}
	By the Taylor expansion for $e^{-i \frac{\tau}{\varepsilon} V}$, we have
	\begin{align}\label{PNStau}
		\begin{aligned}
			P_N\left(\mathcal{S}_\tau\left(P_N \psi\left(t_n\right)\right)\right)= & e^{i\frac{\varepsilon}{2} \tau \Delta} P_N \psi\left(t_n\right)-i \frac{\tau}{\varepsilon} P_N\left(e^{i\frac{\varepsilon}{2} \frac{\tau}{2} \Delta} V e^{i \frac{\varepsilon}{2}\frac{\tau}{2} \Delta} P_N \psi\left(t_n\right)\right) \\
			& -\frac{\tau^2}{\varepsilon^2} P_N\left(\int_0^1(1-\theta) e^{i \frac{\varepsilon}{2}\frac{\tau}{2} \Delta} e^{-i \frac{\tau}{\varepsilon}   V} V^2 e^{i \frac{\varepsilon}{2}\frac{\tau}{2} \Delta} P_N \psi\left(t_n\right) d \theta\right) .
		\end{aligned}
	\end{align}		
	By repeatedly using the Duhamel's principle, one can write
	$$
	\begin{aligned}
		P_N \psi\left(t_{n+1}\right)= & P_N\left(e^{i\frac{\varepsilon}{2} \tau \Delta} \psi\left(t_n\right)\right)-i \frac{1}{\varepsilon} P_N\left(\int_0^\tau e^{i\frac{\varepsilon}{2}(\tau-s) \Delta} V e^{i\frac{\varepsilon}{2} s \Delta} \psi\left(t_n\right) d s\right) \\
		& -\frac{1}{\varepsilon^2} P_N\left(\int_0^\tau \int_0^s e^{i\frac{\varepsilon}{2}(\tau-s) \Delta} V e^{i\frac{\varepsilon}{2}(s-w) \Delta} V \psi\left(t_n+w\right) d w d s\right) .
	\end{aligned}
	$$		
	Recall Assumption \ref{AssumpAB} and apply the Fourier projections, 
	\begin{align}\label{PNpsi}
		\begin{aligned}
			P_N \psi\left(t_{n+1}\right)= & e^{i \frac{\varepsilon}{2}\tau \Delta} P_N \psi\left(t_n\right)-i \frac{1}{\varepsilon} \int_0^\tau P_N\left(e^{i\frac{\varepsilon}{2}(\tau-s) \Delta} V e^{i\frac{\varepsilon}{2} s \Delta} P_N \psi\left(t_n\right)\right) d s \\
			& -\frac{1}{\varepsilon^2} \int_0^\tau \int_0^s P_N\left(e^{i\frac{\varepsilon}{2}(\tau-s) \Delta} V e^{i\frac{\varepsilon}{2}(s-w) \Delta} V P_N \psi\left(t_n+w\right)\right) d w d s-r_h^n,
		\end{aligned}
	\end{align}
	with $\left\|r_h^n(x)\right\|_{L^2} \lesssim \frac{\tau}{\varepsilon} h^{m}$. 
	Combining \eqref{PNStau} and \eqref{PNpsi}, the local truncation error of the TSFP \eqref{TSFP} for the Schr\"odinger equation at time $t_n$ can be written as  \cite{Lubich2008}
	\begin{align}\label{BigErr}
		\begin{aligned}
			 \mathcal{E}^n &:= P_N \mathcal{S}_\tau(P_N \psi(t_n))-P_N \psi(t_{n+1})\\
			& = P_N \mathcal{F}(P_N \psi(t_n))-\frac{\tau^2}{2\varepsilon^2} B^n\left(\frac{\tau}{2}, \frac{\tau}{2}\right)+\frac{r_1^n}{\varepsilon^2}   +\frac{1}{\varepsilon^2} \int_0^\tau\int_0^s B^n(s, w) d w d s+\frac{r_2^n}{\varepsilon^2} +r_h^n,
		\end{aligned}
	\end{align}
	where 
	$$
	\begin{aligned}
		& \mathcal{F}\left(P_N \psi\left(t_n\right)\right)=-i \frac{\tau}{\varepsilon} f^n\left(\frac{\tau}{2}\right)+i \frac{1}{\varepsilon} \int_0^\tau f^n(s) d s, \ f^n(s)=e^{i\frac{\varepsilon}{2}(\tau-s) \Delta} V e^{i \frac{\varepsilon}{2}s \Delta} P_N \psi\left(t_n\right)\\
		& B^n(s, w)=P_N\left(e^{i\frac{\varepsilon}{2}(\tau-s) \Delta} V e^{i\frac{\varepsilon}{2}(s-w) \Delta} V e^{i\frac{\varepsilon}{2} w \Delta} P_N \psi\left(t_n\right)\right),\\
		& r_1^n=-\tau^2 \int_0^1(1-\theta) P_N\left(e^{i \frac{\varepsilon}{2} \frac{\tau}{2} \Delta}\left(e^{-i\theta \frac{\tau}{\varepsilon} V}-1\right) V^2 e^{i\frac{\varepsilon}{2} \frac{\tau}{2} \Delta} P_N \psi\left(t_n\right)\right) d \theta, \\
		& r_2^n=\int_0^\tau \int_0^s\left(P_N\left(e^{i\frac{\varepsilon}{2}(\tau-s) \Delta} V e^{i\frac{\varepsilon}{2}(s-w) \Delta} V P_N \psi\left(t_n+w\right)\right)-B^n(s, w)\right) d w d s.		
	\end{aligned}
	$$
	The first term in \eqref{BigErr} can be estimated by the midpoint quadrature rule as \cite{Lubich2008}
	$$
	\left\|\mathcal{F}\left(P_N \psi\left(t_n\right)\right)\right\|_{L^2} \lesssim \frac{\tau^3}{\varepsilon}\left\|[\Delta,[\Delta, V]] P_N \psi\left(t_n\right)\right\|_{L^2} \lesssim \frac{\tau^3}{\varepsilon}\|V\|_{H^4}\left\|\psi\left(t_n\right)\right\|_{H^2},
	$$
	where $[\Delta,[\Delta, V]]$ is the double commutator. 
	One has
	$$
	\left\|r_1^n\right\|_{L^2} \lesssim  \frac{\tau^3}{\varepsilon}\|V\|_{L^2}^3\left\|\psi\left(t_n\right)\right\|_{L^2} \lesssim \frac{\tau^3}{\varepsilon} .
	$$	
		since $e^{i \frac{\varepsilon}{2}\frac{\tau}{2} \Delta}$ preserves the $L^2$-norm and $\left\|\left(e^{-i \theta \frac{\tau}{\varepsilon} V}-1\right) V^2\right\|_{L^2} \lesssim \frac{\tau}{\varepsilon} \theta\|V\|_{L^2}^3$.
	The following estimates are standard \cite{Lubich2008},
	$$
	\begin{aligned}
		& \left\|r_2^n\right\|_{L^2} \lesssim \frac{\tau^3}{\varepsilon}\|V\|_{L^2}^2\|V \psi(\cdot)\|_{L^{\infty}\left([0, \tau] ; L^2\right)} \lesssim \frac{\tau^3}{\varepsilon}, \\
		& \left\|-\frac{\tau^2}{2} B\left(\frac{\tau}{2}, \frac{\tau}{2}\right)+\int_0^\tau \int_0^s B(s, w) d w d s\right\|_{L^2} \lesssim \tau^3\|V\|_{H^2}^2\left\|\psi\left(t_n\right)\right\|_{H^2} \lesssim \frac{\tau^3}{\e^2} .
	\end{aligned}
	$$
	Combining the above estimates together, the following $L^2$-estimate 
	\begin{align}\label{BigEestimateL2}
		||\mathcal{E}^n||_{L^2}  \lesssim \frac{\tau^3}{\varepsilon^4} + \frac{\tau}{\varepsilon} h^{m}
	\end{align}
	holds for $m \geq 2$.
	
     Finally, combining  \eqref{step1} and \eqref{BigEestimateL2} , one obtains
	$$
	\begin{aligned}
		\left\|e^{n+1}\right\|_{L^2} & \leq\left\|I_N \psi^{n+1}-P_N \mathcal{S}_\tau\left(P_N \psi\left(t_n\right)\right)\right\|_{L^2}+\left\|\mathcal{E}^n\right\|_{L^2} \\
		& \leq\left\|I_N \psi^{(**)}-I_N \psi^{\langle **\rangle}\right\|_{L^2}+\left\|I_N \psi^{\langle **\rangle}-P_N \psi^{\langle **\rangle}\right\|_{L^2} +\left\|\mathcal{E}^n\right\|_{L^2}\\
		& \leq\left\|e^n\right\|_{L^2}+C_4\left(\frac{\tau}{\varepsilon} h^m+\frac{\tau^3}{\varepsilon^4}\right) .
	\end{aligned}
	$$		
	Together with $\left\|e^0\right\|_{L^2} \leq C_2 h^m$, the following estimate 
	\begin{align}\label{INpsi}
		\left\|e^{n+1}\right\|_{L^2} \leq C_4  \frac{t_{n+1}}{\varepsilon} \left(h^m+\frac{\tau^2}{\varepsilon^3}\right)+C_2 h^m, \quad 0 \leq n \leq \frac{T}{\tau}-1,
	\end{align}
	holds. From \eqref{projectionErr} and \eqref{INpsi}, the estimate can be derived		
	$$
	\left\|I_N \psi^n-\psi\left(t_n\right)\right\|_{L^2} 
	\leq C_0 G_m \frac{T}{\varepsilon} \left(h^m+\frac{\tau^2}{\varepsilon^3}\right) .
	$$
	Here $C_0 = \max\{C_1+C_2,C_4\}$  is independent of $\varepsilon, h, \tau, n, m$, and $G_m$ is independent of $\varepsilon, h, \tau$ and only depends on $V$ and $\psi$. This completes the proof of $L^2$-estimate \eqref{H1TSFP}.

\section{Proof of \eqref{ck-norm}} \label{app:B}
By the definition of projection onto $\mathscr{U}_j^L\left(\gamma_K\right)$ for each $\boldsymbol{U}_j^L(z)(j=1,2)$ given in \eqref{Online-eq},
\small{$$
\left(\mathcal{P}_{\mathscr{U}_j^L(\gamma K)}\left[\boldsymbol{U}_j^L(z)\right]\right)^2=\sum_{k, s=1}^K c_k(z) c_s(z) \boldsymbol{U}_j^L\left(z_k\right) \boldsymbol{U}_j^L\left(z_s\right),
$$}
thus
\small{$$\begin{aligned}
	\int_{\Omega}\left(\mathcal{P}_{\mathscr{U}_j^L\left(\gamma_K\right)}\left[\boldsymbol{U}_j^L(z)\right]\right)^2 d x&=\sum_{k, s=1}^K c_k(z) c_s(z) \int_{\Omega} \boldsymbol{U}_j^L\left(z_k\right) \boldsymbol{U}_j^L\left(z_s\right) d x\\
	&:=\mathbf{c}_j^T \mathbf{G}_j^L \mathbf{c}_j \geq \lambda_{j, 0}\left\|\mathbf{c}_j\right\|^2,
\end{aligned}$$}
where $\mathrm{G}_j^L$ is the Gramian matrix of $\mathscr{U}_j^L\left(\gamma_K\right)$ defined in \eqref{Gramian} (where $\langle\cdot, \cdot\rangle$ applies to $L_x^2$ here) and $\lambda_{j, 0}>0$ is its minimum eigenvalue. 
Since for all $z$,
$$
\int_{\Omega}\left(\mathcal{P}_{\mathscr{U}_j^L\left(\gamma_K\right)}\left[\boldsymbol{U}_j^L(z)\right]\right)^2 d x \leq \int_{\Omega}\left[\boldsymbol{U}_j^L(z)\right]^2 d x,\, j=1,2,
$$
then
\small{\begin{align}\label{c-norm}
	\left\|\mathbf{c}_j\right\| \leq \frac{1}{\sqrt{\lambda_{0, j}}}\left(\int_{\Omega}\left(\mathcal{P}_{\mathscr{U}_j^L\left(\gamma_K\right)}\left[\boldsymbol{U}_j^L(z)\right]\right)^2 d x\right)^{1 / 2} \leq \frac{1}{\sqrt{\lambda_{j, 0}}}\left\|\boldsymbol{U}_j^L(z)\right\|_0,
\end{align}}
for all $z$. Thus
$$
\|\| \mathbf{c}_j\|\|_{L_z^2} \leq \frac{1}{\sqrt{\lambda_{j, 0}}}\left\|\boldsymbol{U}_j^L(z)\right\|_{L_x^2 L_z^2} .
$$
Since $\sqrt{x}$ is a monotone function,
$$
\|\| \mathbf{c}_j\|\|_{L_z^2}=\left(\sum_{k=1}^K\left\|c_k(z)\right\|_{L_z^2}^2\right)^{1 / 2} .
$$
By the regularity of $\boldsymbol{U}_j^L$ and the assumption that the volume of $\Omega$ is bounded, \eqref{c-norm} implies \eqref{ck-norm}.

\begin{figure}[!ht]
    \centering
    \includegraphics[width=0.8\textwidth]{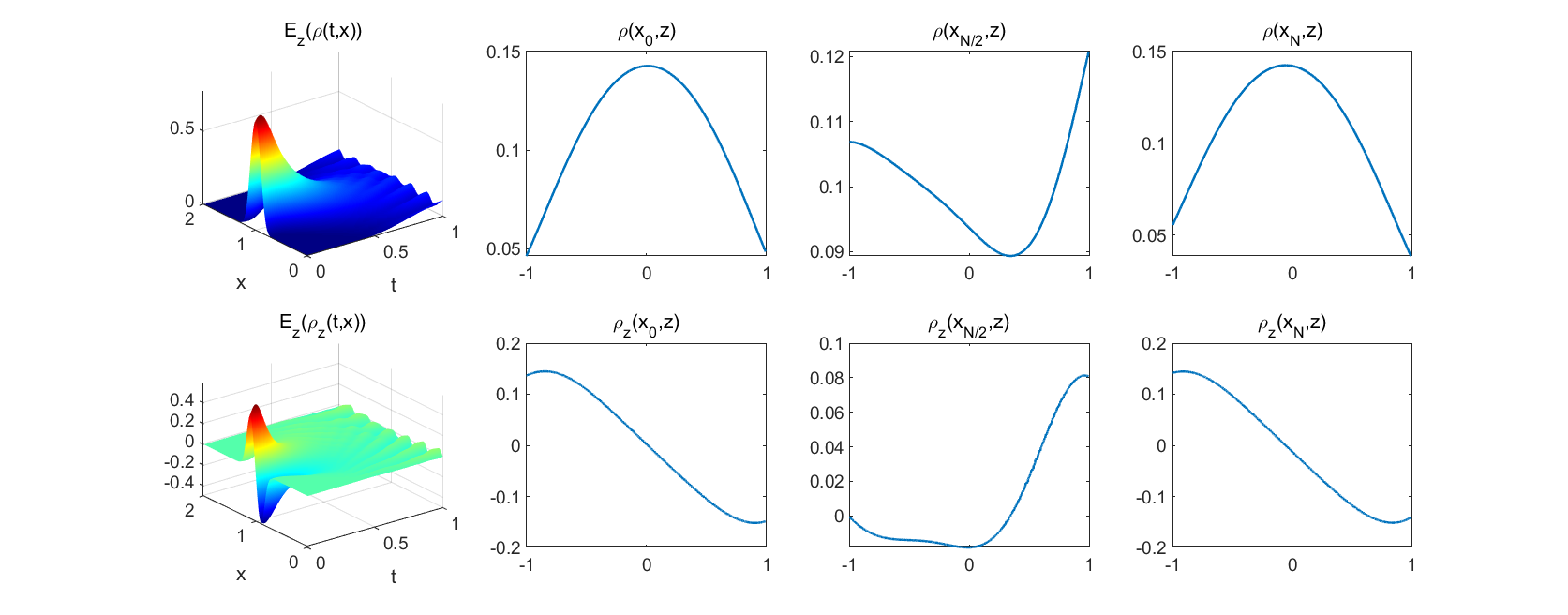}
    \includegraphics[width=0.8\textwidth]{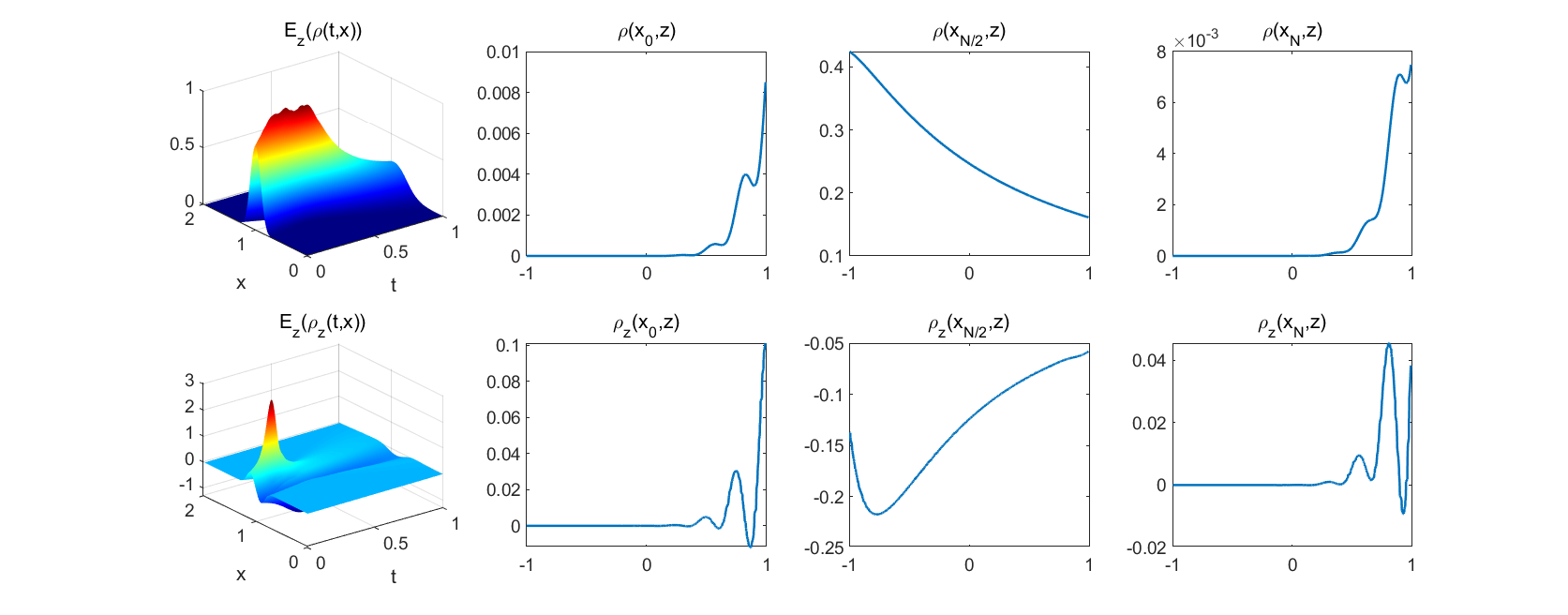}
    \includegraphics[width=0.8\textwidth]{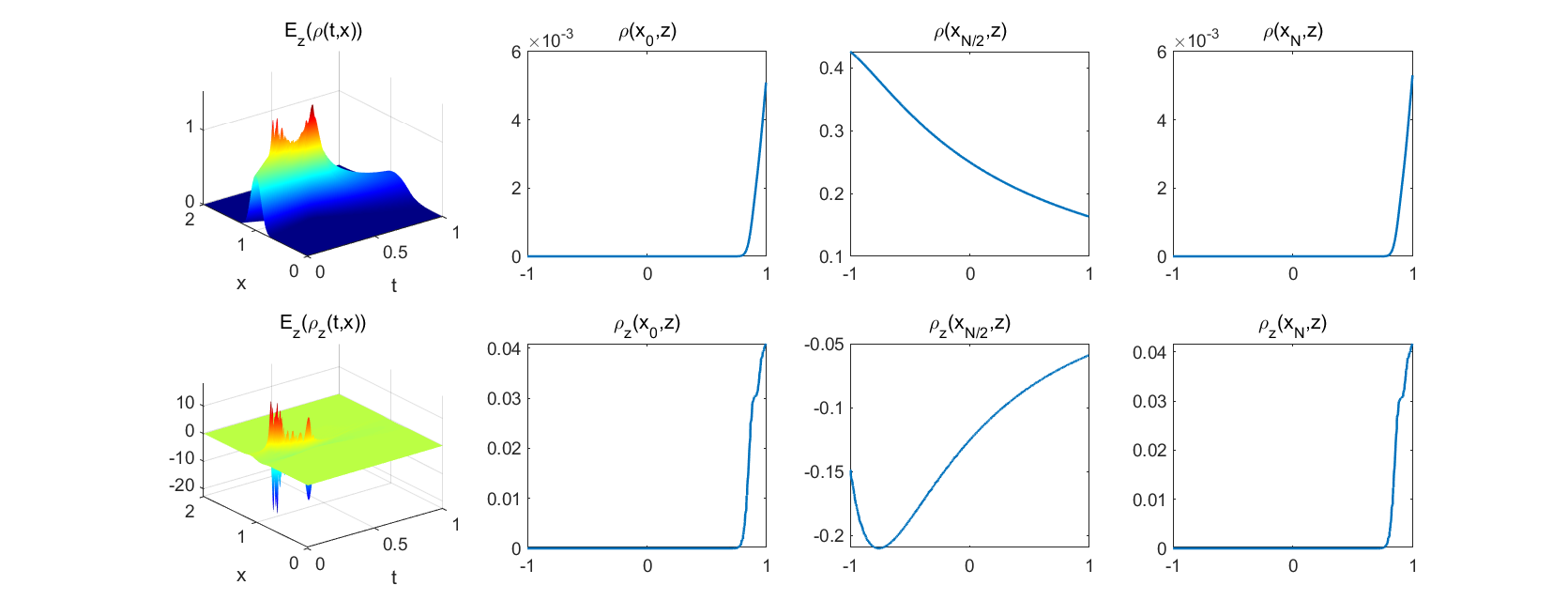}
    \caption{Test 1. The mean of $\rho$ and $\rho_z$ (first column) and the behavior of $\rho$ and $\rho_z$ in the random space at three different points $x_0, x_{N/2}, x_N$ (second to fourth columns) at $T=1$. The first two rows are the results for $\e = 0.1$; the medium two rows are for $\e = 0.01$; and the last two rows are for $\e=0.001$.}
    \label{fig:compare}
\end{figure}

\section{Numerical observation on the bound of $||\rho_z||_{L^2}$}

Regarding Test 1 in our numerical examples, we investigate $\rho$ and $\rho_z$ and study their behaviours in the random space by plotting some quantities in Figure \ref{fig:compare}. Based on these observations, we conclude that the $L_z^2$ norm of $\rho_z$ is bounded by $O(1/\e)$. 

If one considers the WKB ansatz that solves the Schr\"odinger equation \eqref{Schro-eqn}, i.e.,    
\begin{equation}\label{eq:u_wkb}   
 \psi^{\varepsilon}(t,x,z)=\sqrt{\rho^{\varepsilon}(t,x,z)} \exp \left(\frac{i}{\varepsilon} S^{\varepsilon}(t,x,z)\right),
 \end{equation}
where $\rho^{\varepsilon}=\left|\psi^{\varepsilon}\right|^2$ and $S^{\varepsilon}$ are the density and phase function respectively. Inserting \eqref{eq:u_wkb} into the Schr\"odinger equation \eqref{Schro-eqn} and separating the real and imaginary parts give 
 $$\begin{aligned}& \rho_t^{\varepsilon}+\partial_x\left(\rho^{\varepsilon} S_x^{\varepsilon}\right)=0, \\& S_t^{\varepsilon}+\frac{1}{2}\left|S_x^{\varepsilon}\right|^2+V^{\varepsilon}(x,z)=\frac{\varepsilon^2}{2} \frac{\left(\sqrt{\rho^{\varepsilon}}\right)_{xx}}{\sqrt{\rho^{\varepsilon}}} .\end{aligned} 
 $$

Since this coupled system is nonlinear, even just for the WKB form of the wave solution, a rigorous analysis on studying the regularity of $\rho$ in the random space is non-trivial, thus we consider it as a future work and only give findings based on our numerical observations.

\bibliographystyle{siamplain}
\bibliography{references}
\end{document}